\documentclass{amsart}
\usepackage{amsmath}
\usepackage{amssymb}
\usepackage{amsthm}
\usepackage[msc-links,abbrev]{amsrefs}
\usepackage{hyperref}
\usepackage{mathrsfs}
\usepackage{mathtools}
\usepackage{slashed}
\usepackage{tikz-cd}
\usepackage{enumitem}

\setenumerate{label=(\roman*)}

\DeclareMathOperator{\divsymb}{div}

\DeclareMathOperator{\Id}{Id}

\DeclareMathOperator{\im}{im}
\DeclareMathOperator{\tr}{tr}
\DeclareMathOperator{\Vol}{Vol}
\DeclareMathOperator{\dvol}{dV}

\DeclareMathOperator{\darea}{dA}

\DeclareMathOperator{\Ric}{Ric}

\DeclareMathOperator{\Rm}{Rm}

\DeclareMathOperator{\Pf}{Pf}

\DeclareMathOperator{\fp}{fp}
\DeclareMathOperator{\Div}{DIV}
\DeclareMathOperator{\Integral}{INT}
\DeclareMathOperator{\ConfInv}{CONF}

\newcommand{\ambientstraight}{straight}

\newcommand{\straight}{straightenable}
\newcommand{\Straight}{Straightenable}

\newcommand{\Sch}{\mathsf{P}}
\newcommand{\trSch}{\mathsf{J}}
\newcommand{\Weyl}{W}

\newcommand{\defn}[1]{{\boldmath\bfseries#1}}

\newcommand{\Rint}{\sideset{^{R}\!\!}{}\int}
\newcommand{\smallRint}{{}^R\!\!\int}

\newcommand{\cg}{\widetilde{g}}

\newcommand{\cI}{\widetilde{I}}

\newcommand{\cL}{\widetilde{L}}

\newcommand{\cR}{\widetilde{R}}
\newcommand{\cS}{\widetilde{S}}
\newcommand{\cT}{\widetilde{T}}
\newcommand{\cU}{\widetilde{U}}

\newcommand{\cX}{\widetilde{X}}

\newcommand{\cnabla}{\widetilde{\nabla}}

\newcommand{\cDelta}{\widetilde{\Delta}}

\newcommand{\cmE}{\widetilde{\mathcal{E}}}

\newcommand{\cmG}{\widetilde{\mathcal{G}}}

\newcommand{\cmW}{\widetilde{\mathcal{W}}}

\newcommand{\cRic}{\widetilde{\Ric}}
\DeclareMathOperator{\cRm}{\widetilde{\Rm}}

\newcommand{\hu}{\widehat{u}}

\newcommand{\hT}{\widehat{T}}

\newcommand{\hDelta}{\widehat{\Delta}}

\newcommand{\hnabla}{\widehat{\nabla}}

\newcommand{\hpi}{\widehat{\pi}}

\newcommand{\hmG}{\widehat{\mathcal{G}}}

\newcommand{\lp}{\langle}
\newcommand{\rp}{\rangle}
\newcommand{\lv}{\lvert}
\newcommand{\rv}{\rvert}

\newcommand{\mG}{\mathcal{G}}

\newcommand{\mI}{\mathcal{I}}

\newcommand{\mP}{\mathcal{P}}

\newcommand{\mU}{\mathcal{U}}
\newcommand{\mV}{\mathcal{V}}
\newcommand{\mW}{\mathcal{W}}

\newcommand{\bN}{\mathbb{N}}

\newcommand{\bR}{\mathbb{R}}

\newcommand{\sF}{\mathscr{F}}

\newcommand{\onf}{\mathsf{n}}

\def\sideremark#1{\ifvmode\leavevmode\fi\vadjust{\vbox to0pt{\vss
 \hbox to 0pt{\hskip\hsize\hskip1em
 \vbox{\hsize3cm\tiny\raggedright\pretolerance10000
 \noindent #1\hfill}\hss}\vbox to8pt{\vfil}\vss}}}

\newcommand{\suchthat}{\mathrel{}:\mathrel{}}

\newtheorem{theorem}{Theorem}[section]
\newtheorem{proposition}[theorem]{Proposition}
\newtheorem{lemma}[theorem]{Lemma}
\newtheorem{corollary}[theorem]{Corollary}

\theoremstyle{definition}
\newtheorem{definition}[theorem]{Definition}

\newtheorem{example}[theorem]{Example}

\theoremstyle{remark}
\newtheorem{remark}[theorem]{Remark}

\numberwithin{equation}{section}

\usepackage{tikz-cd}
\usetikzlibrary{arrows}

\begin{document}

\title[Computing renormalized curvature integrals]{Computing renormalized curvature integrals on Poincar\'e--Einstein manifolds}
\author[J.\ S.\ Case]{Jeffrey S.\ Case}
\address{Department of Mathematics \\ Penn State University \\ University Park, PA 16802 \\ USA}
\email{jscase@psu.edu}
\author[A. Khaitan]{Ayush Khaitan}
\address{Department of Mathematics \\ Rutgers University \\ Hill Center for the Mathematical Sciences \\ 110 Frelinghuysen Rd. \\ Piscataway, NJ 08854 }
\email{ayush.khaitan@rutgers.edu}
\author[Y.-J.\ Lin]{Yueh-Ju Lin}
\address{Department of Mathematics, Statistics, and Physics \\ Wichita State University \\ Wichita, KS 67260 \\ USA}
\email{yueh-ju.lin@wichita.edu}
\author[A.\ J.\ Tyrrell]{Aaron J.\ Tyrrell}
\address{Department of Mathematics \\ University of Notre Dame \\ Notre Dame, IN 46556 \\ USA}
\email{atyrrell@nd.edu}
\author[W.\ Yuan]{Wei Yuan}
\address{Department of Mathematics \\ Sun Yat-sen University \\ Guangzhou, Guangdong 510275 \\ China}
\email{yuanw9@mail.sysu.edu.cn}
\begin{abstract}
 We describe a general procedure for computing renormalized curvature integrals on Poincar\'e--Einstein manifolds.
 In particular, we explain the connection between the Gauss--Bonnet-type formulas of Albin and Chang--Qing--Yang for the renormalized volume, and explicitly identify a scalar conformal invariant in the latter formula.
 Our approach constructs scalar conformal invariants of weight $-n$ on $n$-manifolds, $n \geq 8$, that are natural divergences;
 these imply that the scalar invariant in the Chang--Qing--Yang formula is not unique in dimension $n \geq 8$.
 Our procedure also produces explicit conformally invariant Gauss--Bonnet-type formulas for compact Einstein manifolds.
\end{abstract}
\maketitle

\section{Introduction}
\label{sec:intro}

Poincar\'e--Einstein manifolds are complete Einstein manifolds with a well-defined conformal boundary;
the model example is the Poincar\'e ball model of hyperbolic space.
The rich interplay between the conformal geometry of the boundary of a Poincar\'e--Einstein manifold and the Riemannian geometry of its interior, which is at the heart of the AdS/CFT correspondence~\cites{Maldacena1998,Witten1998}, leads to many constructions of local and nonlocal conformal invariants (e.g.\ \cites{GrahamZworski2003,GonzalezQing2010,Juhl2013,BaileyEastwoodGraham1994}).
Nevertheless, it remains a fundamental open problem to understand the moduli space of Poincar\'e--Einstein manifolds with a given conformal boundary (cf.\ \cites{Anderson2001b,ChangGe2018,ChangGe2020,ChangGeQing2020}).

A first step in studying this moduli space is to develop a good understanding of global invariants of Poincar\'e--Einstein manifolds.
In even dimensions, renormalized curvature integrals~\cite{Albin2009} give many such invariants:
Let $r$ be a geodesic defining function for an even-dimensional Poincar\'e--Einstein manifold $(M,g_+)$.
A \defn{renormalized curvature integral}
\begin{equation*}
 \Rint I \dvol := \fp \int_{r>\varepsilon} I^{g_+} \dvol_{g_+}
\end{equation*}
is the constant term in the Laurent expansion of the total integral of a scalar Riemannian invariant $I$.
We say that the renormalized curvature integral $\smallRint I \dvol$ is \defn{even} if $I$ is even;
i.e.\ if $I$ is invariant under orientation-reversing diffeomorphisms.
In this paper we restrict our attention to \defn{admissible} renormalized curvature integrals;
i.e.\ renormalized curvature integrals of linear combinations of scalar Riemannian invariants of weight $w \geq -n$ on $n$-manifolds.

There are many renormalized curvature integrals of interest.
The most commonly studied such integral is the \defn{renormalized volume}~\cite{Graham2000}
\begin{equation*}
 \mV := \Rint \dvol ;
\end{equation*}
i.e.\ the renormalized integral of the constant scalar $1$.
Another example is the renormalized integral of the Pfaffian $\Pf$ (of the Riemann curvature tensor).
As a generalization of the Gauss--Bonnet--Chern Theorem~\cite{Chern1955}, Albin showed~\cite{Albin2009} that the Euler characteristic is proportional to this renormalized curvature integral:
\begin{equation}
 \label{eqn:renormalized-pfaffian-formula}
 \Rint \Pf \dvol = (2\pi)^{n/2}\chi(M) .
\end{equation}
These are both even admissible renormalized curvature integrals.

In dimension four, the Poincar\'e--Einstein condition implies that any even admissible renormalized curvature integral is a linear combination of the renormalized volume and the renormalized $L^2$-norm of the Weyl tensor $\Weyl$.
Moreover, $\lv \Weyl \rv^2 \dvol$ is conformally invariant in this dimension, and hence $\int \lv \Weyl \rv^2 \dvol$ is finite.
For example, expanding the Pfaffian in Equation~\eqref{eqn:renormalized-pfaffian-formula} recovers Anderson's formula~\cite{Anderson2001b}
\begin{equation}
 \label{eqn:renormalized-volume-4d}
 4\pi^2\chi(M) = 3\mV + \frac{1}{8}\int_M \lv \Weyl \rv^2 \dvol
\end{equation}
for the renormalized volume.

In higher dimensions there are many additional renormalized curvature integrals.
While Equation~\eqref{eqn:renormalized-pfaffian-formula} still relates the renormalized volume to the Euler characteristic, it now does so through the renormalized integral of a scalar polynomial in the Weyl tensor.
In a different direction, Chang, Qing, and Yang showed~\cite{ChangQingYang2006}, in analogy with Equation~\eqref{eqn:renormalized-volume-4d}, that there is a scalar conformal invariant $\mW_n$ of weight $-n$ such that if $(M^n,g_+)$ is Poincar\'e--Einstein, then
\begin{equation}
 \label{eqn:renormalized-volume-cqy}
 (2\pi)^{n/2}\chi(M) = (-1)^{n/2}(n-1)!!\mV + \int_M \mW_n \dvol ;
\end{equation}
since $\mW_n$ has weight $-n$, the integral $\int \mW_n \dvol$ converges.
An explicit formula for $\mW_n$ is only known in dimension at most six~\cites{Graham2000};
in higher dimensions, the existence of $\mW_n$ follows from Alexakis' classification of conformally invariant integrals of scalar Riemannian invariants~\cite{Alexakis2012}.
For this reason it is unclear whether one can derive one of Equation~\eqref{eqn:renormalized-pfaffian-formula} or~\eqref{eqn:renormalized-volume-cqy} from the other.

The main purpose of this paper is to describe a procedure that is independent of Alexakis' classification and explicitly computes many renormalized curvature integrals.
Applying our procedure to the Pfaffian produces an ambient formula for $\mW_n$ and explains the relationship between Equations~\eqref{eqn:renormalized-pfaffian-formula} and~\eqref{eqn:renormalized-volume-cqy}:

\begin{theorem}
 \label{thm:cqy}
 Let $(M^n,g_+)$ be an even-dimensional Poincar\'e--Einstein manifold.
 Then
 \begin{equation*}
  (2\pi)^{n/2}\chi(M) = (-1)^{n/2}(n-1)!!\mV + \sum_{\ell=2}^{n/2} \left(-2\right)^{\ell-n/2}\frac{(\ell-1)!}{(n/2-1)!}\int_M \mP_{\ell,n} \dvol ,
 \end{equation*}
 where $\mP_{\ell,n} := i^\ast \bigl( \cDelta^{n/2-\ell}\Pf_\ell(\cRm) \bigr)$.
\end{theorem}

Here $\Pf_\ell$ is the Pfaffian-like polynomial
\begin{align}
 \label{eqn:generalized-pfaffian} \Pf_\ell(T) & := 2^{-\ell}(2\ell-1)!!\delta_{b_1 \dotsm b_{2\ell}}^{a_1 \dotsm a_{2\ell}} T_{a_1a_2}{}^{b_1b_2} \dotsm T_{a_{2\ell-1}a_{2\ell}}{}^{b_{2\ell-1}b_{2\ell}} , \\
 \label{eqn:generalized-kronecker} \delta_{b_1 \dotsm b_k}^{a_1 \dotsm a_k} & := \delta_{[b_1}^{[a_1} \dotsm \delta_{b_k]}^{a_k]} ,
\end{align}
on $(2,2)$-tensors $T$, with the convention $\Pf_0(T):=1$.
Thus $\Pf_{n/2}(\Rm)$ is the Pfaffian of an $n$-manifold and $k!\delta_{b_1 \dotsm b_k}^{a_1 \dotsm a_k}$ is the generalized Kronecker delta.
With this notation, $\Pf_\ell(\cRm)$ is the scalar invariant determined by applying $\Pf_\ell$ to the Riemann curvature tensor $\cRm$ of the ambient space of $(M^n,[g_+])$, and $\mP_{\ell,n}$ is the evaluation at $(M^n,g_+)$ of the scalar conformal invariant of weight $-n$ determined by the scalar invariant $\cDelta^{n/2-\ell}\Pf_{\ell}(\cRm)$ built using also the Laplacian $\cDelta$ in the ambient space.
See Section~\ref{sec:bg} for a more detailed discussion.

Direct computation gives $\mP_{2,4} = \frac{1}{8}\lv\Weyl\rv^2$, and hence Theorem~\ref{thm:cqy} recovers Anderson's formula~\eqref{eqn:renormalized-volume-4d}.
Section~\ref{sec:algorithm} contains explicit formulas for $\mP_{2,n}$, $\mP_{3,n}$, and $\mP_{4,n}$.
In particular, we recover Chang, Qing, and Yang's formula~\cite{ChangQingYang2006}*{Theorem~4.3} relating the Euler characteristic, renormalized volume, and the integral of a scalar conformal invariant on Poincar\'e--Einstein six-manifolds (see Corollary~\ref{cor:gbc6}), and we obtain the first such formula for Poincar\'e--Einstein eight-manifolds (see Corollary~\ref{cor:gbc8}).

Our proof of Theorem~\ref{thm:cqy} also produces a formula for the Euler characteristic of a compact even-dimensional Einstein manifold as a linear combination of its volume and the integral of a scalar conformal invariant:

\begin{theorem}
 \label{thm:gbc}
 Let $(M^n,g)$ be a compact even-dimensional Einstein manifold with $\Ric_g = 2\lambda(n-1) g$.
 Then
 \begin{equation*}
  (2\pi)^{n/2}\chi(M) = (2\lambda)^{n/2}(n-1)!!\Vol_g(M) + \sum_{\ell=2}^{n/2} \left( -2 \right)^{\ell-n/2} \frac{(\ell-1)!}{(n/2-1)!}\int_M \mP_{\ell,n} \dvol ,
 \end{equation*}
 where $\mP_{\ell,n}$ is as in Theorem~\ref{thm:cqy}.
\end{theorem}

Theorems~\ref{thm:cqy} and~\ref{thm:gbc} are special cases of a general procedure for computing renormalized curvature integrals.
This procedure begins with a distinguished class of scalar Riemannian invariants that are well-suited to computations in the Fefferman--Graham ambient space:

\begin{definition}
 \label{defn:straight}
 A scalar Riemannian invariant $\cI$ of weight $w$ on ambient $(n+2)$-manifolds is \defn{\ambientstraight} if there is a scalar Riemannian invariant $I$ of weight $w$ on $n$-manifolds such that $\cI\,{}^{\cg} = \tau^w\pi^\ast I^g$ whenever $(\cmG,\cg)$,
 \begin{equation}
  \label{eqn:straight-and-normal}
  \begin{split}
  \cmG & := (0,\infty)_t \times M \times (-\varepsilon,\varepsilon)_\rho , \\
  \cg & := 2\rho \, dt^2 + 2t \, dt \, d\rho + \tau^2 g , \\
  \tau & := t(1+\lambda\rho) ,
  \end{split}
 \end{equation}
 is the canonical straight and normal ambient space associated to an Einstein $n$-manifold $(M,g)$ with $\Ric = 2\lambda(n-1) g$, where $\pi \colon \cmG \to M$ is the canonical projection.
 In this case we call $I$ a \defn{\straight} invariant associated to $\cI$.
\end{definition}

For example, computations of Fefferman and Graham~\cite{FeffermanGraham2012} imply that $\Pf_\ell(W)$ is a \straight\ invariant associated to the \ambientstraight\ invariant $\Pf_\ell(\cRm)$.
See Section~\ref{sec:bg} for a discussion of scalar (and tensor) Riemannian invariants.
The sense in which $(\cmG,\cg)$ is canonical is explained by Fefferman and Graham~\cite{FeffermanGraham2012}.

There are two key properties of \ambientstraight\ and \straight\ scalar Riemannian invariants used in this paper.
First, if $I$ is a \straight\ scalar Riemannian invariant associated to the \ambientstraight\ scalar Riemannian invariant $\cI$, then the evaluations of $I$ and $i^\ast\cI$ at an Einstein manifold agree.
Second, there are two recursive constructions of \ambientstraight\ invariants via derivatives;
see Section~\ref{sec:conformal-invariant}.
These properties allow us to construct complicated scalar conformal invariants (e.g.\ $\mP_{\ell,n}$) that are nevertheless easy to evaluate for Einstein metrics.

Remarkably, our second recursive construction produces, in every even dimension $n \geq 8$, nontrivial scalar conformal invariants of weight $-n$ on $n$-manifolds that are natural divergences;
see Proposition~\ref{prop:conf-inv-div} for a precise statement and Equation~\eqref{eqn:dim8-examples} for explicit examples in dimension $n=8$.
Chen and Lu gave~\cite{ChenLu2024} an alternative physical explanation for the existence of these invariants in dimension $n=8$.
Our invariants show that Alexakis' decomposition~\cite{Alexakis2012}
\begin{equation*}
 \Integral = \left\lp\Pf\right\rp \oplus \left( \ConfInv + \Div \right)
\end{equation*}
of the space $\Integral$ of scalar Riemannian invariants whose integrals are conformally invariant on compact $n$-manifolds into the span $\left\lp\Pf\right\rp$ of the Pfaffian, the space $\ConfInv$ of scalar conformal invariants of weight $-n$, and the space $\Div$ of natural divergences is \emph{not} as a direct sum when $n \geq 8$.
In particular, the invariant $\mW_n$ in Equation~\eqref{eqn:renormalized-volume-cqy} is \emph{not} uniquely determined in dimension $n \geq 8$.

The general result underlying Theorems~\ref{thm:cqy} and~\ref{thm:gbc} is the following formula for the renormalized integral of a \straight\ scalar Riemannian invariant in terms of the \emph{convergent} integral of a scalar conformal invariant:

\begin{theorem}
 \label{thm:main-thm}
 Let $n \in \bN$ be an even integer and let $I$ be a \straight\ scalar Riemannian invariant of weight $-2k \in [ -n , 0)$.
 If $(M^n,g_+)$ is a Poincar\'e--Einstein manifold, then
 \begin{equation*}
  \Rint I \dvol = \frac{2^{k-n/2}(k-1)!}{(n/2-1)!(n-2k-1)!!} \int_M \mI_{n/2-k} \dvol ,
 \end{equation*}
 where $\mI_{n/2-k} := i^\ast( \cDelta^{n/2-k}\cI)$ and $\cI$ is a \ambientstraight\ scalar Riemannian invariant to which $I$ is associated.
\end{theorem}

Applying Theorem~\ref{thm:main-thm} and its proof to $\Pf_\ell(\Weyl)$ recovers Theorems~\ref{thm:cqy} and~\ref{thm:gbc}.

We emphasize that $\mI_{n/2-k}$ is a scalar conformal invariant of weight $-n$, and hence $\int \mI_{n/2-k}\dvol$ converges.
Indeed, Theorem~\ref{thm:main-thm} gives a new proof (cf.\ \cite{Albin2009}) that $\smallRint I \dvol$ is independent of the choice of geodesic defining function when $I$ is \straight.
As noted above, there are many \straight\ scalar Riemannian invariants, and so Theorem~\ref{thm:main-thm} is widely applicable;
see Section~\ref{sec:algorithm} for examples.

There are two key observations in the proof of Theorem~\ref{thm:main-thm}.
First, if $\cI$ is a \ambientstraight\ scalar Riemannian invariant, then so too is $\cDelta\cI$;
this is the first recursive construction mentioned above.
Moreover, if $I$ is associated to $\cI$, then $\Delta I + cRI$ is associated to $\cDelta \cI$ for a constant $c$ depending only on the dimension and the weight of $I$;
here $R$ is the scalar curvature.
Iteratively applying this to a \straight\ scalar Riemannian invariant $I$ of weight $-2k \geq -n$ produces a \straight\ scalar Riemannian invariant of weight $-n$ which evaluates to a linear combination of $I$ and divergences on Einstein manifolds.
Second, a parity argument of Albin~\cites{Albin2009} can be adapted to prove that the renormalized integral of a natural divergence is zero.

This paper is organized as follows:

In Section~\ref{sec:bg} we fix our conventions and recall some necessary facts about the Fefferman--Graham ambient space.

In Section~\ref{sec:conformal-invariant} we produce many examples of \ambientstraight\ scalar Riemannian invariants.
As an application, we construct nontrivial scalar conformal invariants of weight $-n$ on $n$-manifolds that are natural divergences.

In Section~\ref{sec:cvi} we prove Theorem~\ref{thm:main-thm}.

In Section~\ref{sec:algorithm} we prove Theorems~\ref{thm:cqy} and~\ref{thm:gbc}, derive more explicit formulas for the invariants $\mP_{\ell,n}$ when $\ell \leq 4$, and discuss how to apply Theorem~\ref{thm:main-thm} to compute more general renormalized curvature integrals.

\section{Conventions and background}
\label{sec:bg}

Let $M^n$ be a smooth manifold of dimension $n \geq 3$.
We typically perform computations using abstract index notation, using lowercase Latin letters to denote factors of $T^\ast M$ (when subscripts) or its dual $TM$ (when superscripts), with repeated indices denoting a contraction via the canonical pairing of $TM$ and $T^\ast M$.
For example, we often write $T_{abc}$ to denote a section $T \in \Gamma(\otimes^3 T^\ast M)$ and $X^a$ to denote a vector field $X \in \Gamma(TM)$.
If also $Y,Z \in \Gamma(TM)$, then
\begin{equation*}
 T(X,Y,Z) = X^a Y^b Z^c T_{abc} .
\end{equation*}
We use square brackets and round parentheses to denote skew-symmetrization and symmetrization, respectively;
e.g.
\begin{align*}
 T_{[abc]} & := \frac{1}{6}(T_{abc} - T_{acb} + T_{bca} - T_{bac} + T_{cab} - T_{cba}) , \\
 T_{(abc)} & := \frac{1}{6}(T_{abc} + T_{acb} + T_{bca} + T_{bac} + T_{cab} + T_{cba}) .
\end{align*}

Let $g$ be a pseudo-Riemannian metric on $M$.
Denote by $g^{-1}$ or $g^{ab}$ the inverse of $g$.
We denote by $\Rm$ or $R_{abcd}$ the Riemann curvature tensor, defined by the convention
\begin{equation*}
 \nabla_a \nabla_b \tau_c - \nabla_b \nabla_a \tau_c = R_{abc}{}^d\tau_d .
\end{equation*}
We use $g_{ab}$ and $g^{ab}$ to lower and raise indices;
e.g.\ $R_{abcd} = R_{abc}{}^eg_{ed}$.
The \defn{Schouten tensor} is
\begin{equation*}
 \Sch_{ab} := \frac{1}{n-2}\left( R_{ab} - \trSch g_{ab} \right) ,
\end{equation*}
where $R_{ab} := R_{acb}{}^c$ is the Ricci tensor, $R := R_a{}^a$ is the scalar curvature, and $\trSch := \Sch_a{}^a$ is the trace of the Schouten tensor.
Recall that the \defn{Weyl tensor}
\begin{equation*}
 \Weyl_{abcd} := R_{abcd} - 2\Sch_{a[c}g_{d]b} + 2\Sch_{b[c}g_{d]a}
\end{equation*}
is conformally covariant: $\Weyl^{e^{2u}g} = e^{2u}\Weyl^g$.
The \defn{Laplacian} is $\Delta := \nabla^a\nabla_a$.

A (even) \defn{tensor Riemannian invariant} $I$ of rank $k$ on $n$-manifolds is an assignment to each pseudo-Riemannian $n$-manifold $(M,g)$ of a section $I^g \in \Gamma(\otimes^kT^\ast M)$ that can be expressed as an $\bR$-linear combination of partial contractions of
\begin{equation}
 \label{eqn:weyl-form}
 \nabla^{N_1}\Rm \otimes \dotsm \otimes \nabla^{N_j}\Rm \otimes \mathop{g^{\otimes\ell}} ,
\end{equation}
where we regard $\nabla^N\Rm \in \Gamma(\otimes^{N+4}T^\ast M)$ and $g \in \Gamma(\otimes^2T^\ast M)$, so that all contractions are made using the inverse metric $g^{-1}$, and by convention the empty tensor product is the constant scalar $1$.
We say that $I$ has \defn{weight} $w \in \bR$ if $I^{c^2g} = c^wI^g$ for all pseudo-Riemannian $n$-manifolds $(M,g)$ and all constants $c>0$;
in this case we say that $I$ is \defn{homogeneous} and we define its \defn{tensor weight} to be $w-k$.
A (even) \defn{scalar Riemannian invariant} is a tensor Riemannian invariant of rank $0$.
A \defn{scalar conformal invariant} of weight $w \in \bR$ is a scalar Riemannian invariant $I$ such that $I^{u^2g} = u^wI^g$ for all positive $u \in C^\infty(M)$.
A \defn{natural divergence} is the divergence of a natural one-form.

The \defn{straight and normal ambient space} $(\cmG,\cg)$ associated to $(M^n,g)$ is the smooth $(n+2)$-manifold $\cmG := (0,\infty)_t \times M \times (-\varepsilon,\varepsilon)_\rho$ together with the pseudo-Riemannian metric
\begin{equation}
 \label{eqn:ambient-metric}
 \cg = 2\rho \, dt^2 + 2t \, dt \, d\rho + t^2 g_\rho ,
\end{equation}
where $g_\rho$ is a one-parameter family of metrics on $M$ such that $g_0=g$ and
\begin{enumerate}
 \item if $n$ is odd, then $\cRic = O(\rho^\infty)$ along $\{ \rho=0 \}$; while
 \item if $n$ is even, then $\cRic = O(\rho^{n/2-1})$ and $\cR = O(\rho^{n/2})$ along $\{ \rho=0 \}$.
\end{enumerate}
Here $\cT = O(\rho^k)$ along $\{ \rho=0 \}$ if $\left.\partial_\rho^\ell\right|_{\rho=0}\cT=0$ for all $\ell < k$.
Note that $\delta_s^\ast\cg = s^2\cg$ for $\delta_s \colon \cmG \to \cmG$, $s > 0$, the one-parameter family of dilations $\delta_s(t,x,\rho) := (st,x,\rho)$.
Fefferman and Graham proved~\cite{FeffermanGraham2012}*{Theorem~2.9} the existence and uniqueness of the straight and normal ambient space.
Uniqueness means that if $n$ is odd, then $g_\rho$ is uniquely determined modulo $O(\rho^\infty)$; while if $n$ is even, then $g_\rho$ and $\tr_g g_\rho$ are uniquely determined modulo $O(\rho^{n/2})$ and $O(\rho^{n/2+1})$, respectively.
Denote by $i_g \colon M \to \cmG$ the inclusion $i_g(x) := (1,x,0)$, and by $\pi \colon \cmG \to M$ the projection $\pi(t,x,\rho) := x$.
Fefferman and Graham also showed~\cite{FeffermanGraham2012}*{Theorem~2.3} that if $(\cmG_j,\cg_j)$, $j \in \{ 1, 2 \}$ are straight and normal ambient metrics for conformally related metrics $g_2 = e^{2u}g_1 \in [g]$, then, after shrinking $\cmG_j$ if necessary, there is a $\delta_s$-equivariant diffeomorphism $\Phi \colon \cmG_1 \to \cmG_2$ such that $\Phi(t,x,0)=(e^{-u(x)}t,x,0)$ and $\Phi^\ast \cg_2 \equiv \cg_1$ modulo the undetermined terms in the expansions of $g_\rho$.

One application of the straight and normal ambient space is to the construction of scalar conformal invariants:
Let $\cI$ be a scalar Riemannian invariant of weight $w \geq -n$ on ambient $(n+2)$-manifolds;
i.e.\ $\cI$ is a scalar Riemannian invariant that we will only evaluate on ambient spaces of $n$-manifolds.
Then $\mI := i^\ast\cI$ defines a scalar conformal invariant of weight $w$ on $n$-manifolds.
We say that $\mI$ is \defn{determined ambiently} by $\cI$.
All even scalar conformal invariants of weight $w \geq -n$ arise in this way~\cite{BaileyEastwoodGraham1994}.

When computing with the ambient metric, we use capital Latin letters to denote factors of $T\cmG$ and its dual;
and we use the symbols $0$, $a$, and $\infty$ to denote projections onto $T(0,\infty)_t$, $TM$, and $T(-\varepsilon,\varepsilon)_\rho$, respectively, and their duals.
For example, Fefferman and Graham observed~\cite{FeffermanGraham2012}*{Equation~(6.1)} that the nonvanishing projections of the Riemann curvature tensor $\cRm$ of $\cg$ can be computed from the formulas
\begin{equation}
 \label{eqn:ambient-curvature-general}
 \begin{aligned}
  \cR_{abcd} & = t^2 \left( R_{abcd} - g_{a[c}g_{d]b}^\prime + g_{b[c}g_{d]a}^\prime + \rho g_{a[c}^\prime g_{d]b}^\prime \right) , \\
  \cR_{\infty abc} & = t^2\nabla_{[c}g_{b]a}^\prime , \\
  \cR_{\infty ab \infty} & = \frac{t^2}{2}\left( g_{ab}^{\prime\prime} - \frac{1}{2}g^{cd}g_{ac}^\prime g_{bd}^\prime \right) ,
 \end{aligned}
\end{equation}
where ${}^\prime := \partial_\rho$ and the geometric invariants on the right-hand side are computed with respect to $g_\rho$.

\section{\Straight\ scalar Riemannian invariants}
\label{sec:conformal-invariant}

In this section we construct many \straight\ scalar Riemannian invariants.
The main point is that there is an explicit straight and normal ambient space associated to an Einstein manifold which is in fact an \emph{exact} solution of $\cRic = 0$.

\begin{lemma}[\cite{FeffermanGraham2012}*{p.\ 67}]
 \label{lem:fg-einstein}
 The straight and normal ambient space $(\cmG,\cg)$ associated to an Einstein manifold $(M^n,g)$, $n \geq 3$, with $\Ric_g = 2\lambda(n-1) g$ is given by Equations~\eqref{eqn:straight-and-normal}.
 Moreover, $\cRic = 0$.
\end{lemma}

Thus \ambientstraight\ scalar Riemannian invariants are characterized by their evaluations at straight and normal ambient spaces associated to an Einstein manifold.
The ambient space $(\cmG,\cg)$ of Lemma~\ref{lem:fg-einstein} is canonical in the sense that if $g_i$, $i \in \{ 1, 2 \}$, are conformally equivalent Einstein metrics on $M$ and if $(\cmG_i,\cg_i)$ are the corresponding ambient spaces from Lemma~\ref{lem:fg-einstein}, then there is a diffeomorphism $\Phi \colon \cmG_1 \to \cmG_2$ that fixes $\mG := \{ \rho = 0 \}$ and for which $\cg_1 - \Phi^\ast \cg_2$ vanishes to infinite order along $\mG$~\cite{FeffermanGraham2012}*{Proposition~7.5}.

A \straight\ scalar Riemannian invariant need not be conformally invariant;
e.g.\ $\lvert \nabla R \rvert^2$ is a \straight\ scalar invariant associated to $0$, but it is clearly not conformally invariant.
Nevertheless, a \straight\ scalar Riemannian invariant is necessarily equal to a conformal invariant when evaluated at an Einstein manifold.

\begin{lemma}
 \label{lem:straight-is-almost-conformal}
 Fix $n \in \bN$.
 Let $I$ be a \straight\ scalar Riemannian invariant of weight $w \geq -n$ on $n$-manifolds that is associated to the \ambientstraight\ scalar Riemannian invariant $\cI$.
 Let $\mI$ be the scalar conformal invariant determined ambiently by~$\cI$.
 If $(M,g)$ is an Einstein $n$-manifold, then $I^g = \mI^g$.
\end{lemma}

\begin{proof}
 Let $(\cmG,\cg)$ be the straight and normal ambient space associated to the Einstein $n$-manifold $(M,g)$.
 Since $i_g^\ast\tau=1$ and $\pi_g \circ i_g = \Id$, we compute that
 \begin{equation*}
  \mI^g = i_g^\ast\cI{}\,^{\cg} = i_g^\ast( \tau^w \pi_g^\ast I^g ) = I^g . \qedhere
 \end{equation*}
\end{proof}

For our general construction of \straight\ scalar Riemannian invariants, it is convenient to extend Definition~\ref{defn:straight} to tensor Riemannian invariants.

\begin{definition}
 \label{defn:straight-tensor}
 A tensor Riemannian invariant $\cT$ of rank $k$ and weight $w$ on ambient $(n+2)$-manifolds is \defn{\ambientstraight} if there is tensor Riemannian invariant $T$ of rank $k$ and weight $w$ on Riemannian $n$-manifolds such that
 \begin{equation*}
  \cT^{\cg} = \tau^w\pi^\ast T^g
 \end{equation*}
 whenever $(\cmG,\cg)$ is the straight and normal ambient space~\eqref{eqn:straight-and-normal} associated to an Einstein $n$-manifold $(M,g)$.
 In this case we say that $T$ is a \defn{\straight} tensor Riemannian invariant associated to $\cT$.
\end{definition}

\Straight\ tensor Riemannian invariants exist:

\begin{lemma}
 \label{lem:weyl-is-straight}
 The Weyl tensor is a \straight\ tensor Riemannian invariant of rank $4$ and weight $2$ associated to the \ambientstraight\ tensor Riemannian invariant $\cRm$.
\end{lemma}

\begin{proof}
 Let $(\cmG,\cg)$ be the straight and normal ambient space of an Einstein $n$-manifold $(M,g)$.
 Equation~\eqref{eqn:ambient-curvature-general} and Lemma~\ref{lem:fg-einstein} imply that $\cRm = \tau^2 \pi_g^\ast W^g$ as sections of $\otimes^4T^\ast\cmG$.
\end{proof}

Partial contraction of tensor products yields additional \straight\ tensor Riemannian invariants.
In what follows, we say that partial contractions of tensors $S_{a_1 \dotsm a_k}$ and $\cS_{A_1 \dotsm A_k}$ of rank $k$ are the \defn{same} if they are obtained by contracting the same pairs of indices and listing the free indices in the same order;
e.g.\ $S_{acde}T_b{}^{cde}$ and $\cS_{ACDE}\cT_B{}^{CDE}$ are the same contractions of $S \otimes T$ and $\cS \otimes \cT$, respectively.

\begin{lemma}
 \label{lem:partial-contraction}
 Suppose that $S$ and $T$ are \straight\ tensor Riemannian invariants of tensor weight $s$ and $t$, respectively, associated to \ambientstraight\ tensor Riemannian invariants $\cS$ and $\cT$, respectively.
 Then any partial contraction of $S \otimes T$ is \straight, of tensor weight $s+t$, and associated to the same partial contraction of $\cS \otimes \cT$.
\end{lemma}

\begin{proof}
 It follows immediately from Definition~\ref{defn:straight-tensor} that $S \otimes T$ is a \straight\ tensor Riemannian invariant of tensor weight $s+t$ associated to $\cS \otimes \cT$.
 Observe that the projection $\cg^{ab} = \tau^{-2}g^{ab}$ when $g$ is Einstein.
 It follows that any partial of $S \otimes T$ is a \straight\ tensor Riemannian invariant of tensor weight $s+t$ that is associated to the same partial contraction of $\cS \otimes \cT$.
\end{proof}

In particular, scalar polynomials in the Weyl tensor are \straight:

\begin{corollary}
 \label{cor:scalar-polynomial}
 Any complete contraction of $W^{\otimes k}$ is a \straight\ scalar Riemannian invariant of weight $-2k$ associated to the same complete contraction of $\cRm^{\otimes k}$.
\end{corollary}

\begin{proof}
 Apply Lemmas~\ref{lem:weyl-is-straight} and~\ref{lem:partial-contraction}.
\end{proof}

We now give two recursive constructions of \ambientstraight\ tensor invariants.

First, a nice property of the ambient Laplacian on functions implies that the Laplacian of a \ambientstraight\ scalar invariant is \ambientstraight:

\begin{proposition}
 \label{prop:easy-straight-derivatives}
 Let $\cI$ be a \ambientstraight\ scalar Riemannian invariant of weight $-2k$ on ambient $(n+2)$-manifolds and let $\ell \in \bN$.
 Then $\cDelta^\ell\cI$ is a \ambientstraight\ scalar Riemannian invariant of weight $-2k-2\ell$.
 Moreover, if $I$ is associated to $\cI$, then
 \begin{equation}
  \label{eqn:defn-Iell}
  I_\ell := \left( \prod_{j=0}^{\ell-1} \left( \Delta - \frac{4(k+j)(n-2k-2j-1)}{n}\trSch \right) \right) I
 \end{equation}
 is associated to $\cDelta^\ell\cI$.
 Additionally, if $(M^n,g)$ is Einstein, then
 \begin{equation}
  \label{eqn:mI-constant-term}
  I_\ell \equiv \left(-\frac{4}{n}\trSch\right)^{\ell}\frac{\bigl(k+\ell-1\bigr)!\bigl(n-2k-1\bigr)!!}{\bigl(k-1\bigr)!\bigl(n-2k-2\ell-1\bigr)!!}I \mod \Delta I, \dotsc, \Delta^{\ell}I .
 \end{equation}
\end{proposition}

\begin{proof}
 Let $(\cmG,\cg)$ be the straight and normal ambient space associated to an Einstein $n$-manifold $(M,g)$ with $\Ric = 2\lambda(n-1) g$.
 Direct computation shows~\cite{CaseLinYuan2022or}*{Lemma~5.1} that if $u \in C^\infty(M)$ and $w \in \bR$, then
 \begin{equation}
  \label{eqn:ambient-laplacian}
  \cDelta (\tau^w \pi^\ast u) = \tau^{w-2}\pi^\ast \bigl( \Delta + 2\lambda w (n+w-1) \bigr)u .
 \end{equation}
 Repeated application of Equation~\eqref{eqn:ambient-laplacian} implies that
 \begin{equation*}
  \cDelta^\ell\cI = \tau^{w-2\ell}\pi^\ast I_\ell ,
 \end{equation*}
 establishing the first two claims.
 For the third claim, we use the fact that $\trSch$ is constant at Einstein manifolds to compute that
 \begin{align*}
  I_\ell & \equiv \left(-\frac{4}{n}\trSch\right)^\ell \left( \prod_{j=0}^{\ell-1} (k+j) \right)\left(\prod_{j=0}^{\ell-1} (n-2k-2j-1) \right)I \mod \Delta I, \dotsc, \Delta^{\ell}I \\
  & \equiv \left( -\frac{4}{n}\trSch \right)^\ell \frac{(k+\ell-1)!(n-2k-1)!!}{(k-1)!(n-2k-2\ell-1)!!}I \mod \Delta I , \dotsc, \Delta^\ell I . \qedhere
 \end{align*}
\end{proof}

Second, a nice property of the ambient divergence on symmetric tensors allows us to produce higher-order \ambientstraight\ tensor Riemannian invariants:

\begin{proposition}
 \label{prop:preserve-divergence}
 Let $\cT_{A_1 \dotsm A_k}$ be a \ambientstraight\ symmetric tensor Riemannian invariant of rank $k$ and weight $w$.
 Then
 \begin{equation*}
  \cU_{A_1 \dotsm A_{k-1}} := \cnabla^B \cT_{A_1 \dotsm A_{k-1} B} + \frac{k-1}{w-2k+2}\cnabla_{(A_1}\cT_{A_2 \dotsm A_{k-1})B}{}^B
 \end{equation*}
 is a \ambientstraight\ symmetric tensor Riemannian invariant of rank $k-1$ and weight $w-2$.
 Moreover, if $T_{a_1 \dotsm a_k}$ is associated to $\cT_{A_1 \dotsm A_k}$, then
 \begin{equation*}
  U_{a_1 \dotsm a_{k-1}} := \nabla^b T_{a_1 \dotsm a_{k-1}b} + \frac{k-1}{w-2k+2}\nabla_{(a_1}T_{a_2 \dotsm a_{k-1})b}{}^b
 \end{equation*}
 is associated to $\cU_{A_1 \dotsm A_{k-1}}$.
\end{proposition}

\begin{proof}
 Clearly $U_{a_1 \dotsm a_{k-1}}$ and $\cU_{A_1 \dotsm A_{k-1}}$ are symmetric.
 
 Let $(M^n,g)$ be such that $\Ric = 2\lambda(n-1) g$ and let $(\cmG,\cg)$ be its straight and normal ambient space.
 Denote $\sigma := 1 + \lambda \rho$.
 Specializing the formulas~\cite{FeffermanGraham2012}*{Equation~(3.16)} for the Christoffel symbols of a metric of the form~\eqref{eqn:ambient-metric} to $\cg$ yields
 \begin{equation*}
  \begin{split}
   \Gamma_{AB}^0 & = \begin{pmatrix} 0 & 0 & 0 \\ 0 & -\lambda\tau g_{ab} & 0 \\ 0 & 0 & 0 \end{pmatrix} , \\
   \Gamma_{AB}^c & = \begin{pmatrix} 0 & t^{-1}\delta_b^c & 0 \\ t^{-1}\delta_a^c & \Gamma_{ab}^c & \lambda\sigma^{-1}\delta^c_a \\ 0 & \lambda\sigma^{-1}\delta^c_b & 0 \end{pmatrix} , \\
   \Gamma_{AB}^\infty & = \begin{pmatrix} 0 & 0 & t^{-1} \\ 0 & \sigma(\lambda\rho - 1)g_{ab} & 0 \\ t^{-1} & 0 & 0 \end{pmatrix} .
  \end{split}
 \end{equation*}
 Combining this with the assumption $\cT = \tau^w \pi^\ast T$ yields
 \begin{align*}
  \cnabla^B \cT_{a_1 \dotsm a_{k-1} B} & = \tau^{w-2}\nabla^b T_{a_1 \dotsm a_{k-1}b} , \\
  \cnabla^B \cT_{0 a_2 \dotsm a_{k-1}B} & = -\frac{\tau^{w-2}}{t} T_{a_2 \dotsm a_{k-1}b}{}^b , \\
  \cnabla^B \cT_{\infty a_2 \dotsm a_{k-1}B} & = -\frac{\lambda\tau^{w-2}}{\sigma}T_{a_2 \dotsm a_{k-1}b}{}^b ,
 \end{align*}
 and all other components vanish.
 We also compute that
 \begin{align*}
  \cnabla_{(a_1}\cT_{a_2 \dotsm a_{k-1})B}{}^B & = \tau^{w-2}\nabla_{(a_{1}} T_{a_2 \dotsm a_{k-1})b}{}^b , \\
  \cnabla_{(0}\cT_{a_2 \dotsm a_{k-1})B}{}^B & = \frac{w-2k+2}{k-1}\frac{\tau^{w-2}}{t}T_{a_2 \dotsm a_{k-1}b}{}^b , \\
  \cnabla_{(\infty}\cT_{a_2 \dotsm a_{k-1})B}{}^B & = \frac{w-2k+2}{k-1}\frac{\lambda\tau^{w-2}}{\sigma}T_{a_2 \dotsm a_{k-1}b}{}^b ,
 \end{align*}
 and all other components vanish.
 The conclusion readily follows.
\end{proof}

Iteratively applying Proposition~\ref{prop:preserve-divergence} to \ambientstraight\ symmetric tensor Riemannian invariants produces nontrivial scalar conformal invariants which are natural divergences on Einstein manifolds.
For example, if $\cT_{AB}$ is a \ambientstraight\ symmetric tensor Riemannian invariant of rank $2$ and weight $w$, then Proposition~\ref{prop:preserve-divergence} implies that
\begin{equation*}
 \cU := \cnabla^A\cnabla^B\cT_{AB} + \frac{1}{w-2}\cnabla^A\cnabla_A\cT_B{}^B
\end{equation*}
is a \ambientstraight\ scalar Riemannian invariant of weight $w-4$.
Moreover, if $T_{ab}$ is associated to $\cT_{AB}$, then
\begin{equation*}
 U := \nabla^a\nabla^b T_{ab} + \frac{1}{w-2}\nabla^a\nabla_aT_b{}^b 
\end{equation*}
is associated to $\cU$.
The largest weight for which $\cU \not= 0$ is $w=-4$:
The vanishing of the ambient Ricci tensor implies that every \ambientstraight\ symmetric tensor Riemannian invariant of rank $2$ and weight $0$ is trivial.
While $\cR_{ACDE}\cR_B{}^{CDE}$ is a nontrivial \ambientstraight\ symmetric tensor Riemannian invariant of rank $2$ and weight $-2$, the second Bianchi identity and the vanishing of the ambient Ricci tensor imply that
\begin{equation*}
 \cnabla^A\cnabla^B( \cR_{ACDE}\cR_B{}^{CDE} ) - \frac{1}{4}\cnabla^A\cnabla_A( \cR_{BCDE}\cR^{BCDE} ) = 0 .
\end{equation*}
There are two \ambientstraight\ symmetric tensor Riemannian invariants of rank $2$ and weight $-4$.
Applying Proposition~\ref{prop:preserve-divergence} twice to these yields the \ambientstraight\ scalar Riemannian invariants
\begin{align*}
 \cU^{(1)} & := \cnabla^A\cnabla^B ( \cR_{ACBD}\cR^{CEFG}\cR^D{}_{EFG} ) , \\
 \cU^{(2)} & := \cnabla^A\cnabla^B ( \cR_{ACDE}\cR_B{}^C{}_{FG}\cR^{DEFG} ) - \frac{1}{6}\cDelta( \cR_{AB}{}^{CD}\cR_{CD}{}^{EF}\cR_{EF}{}^{AB} ) ,
\end{align*}
of weight $-8$.

In fact, direct computation using Equation~\eqref{eqn:ambient-curvature-general} and the formulas~\cite{FeffermanGraham2012}*{Equations~(3.6) and~(3.16)} for $g_\rho$ mod $O(\rho^2)$ and the Christoffel symbols of an ambient metric in straight and normal form shows that the scalar conformal invariants $\mU^{(j)} := i^\ast\cU^{(j)}$, $j \in \{ 1, 2 \}$, of weight $-8$ on $8$-manifolds determined by $\cU^{(j)}$ are
\begin{equation}
 \label{eqn:dim8-examples}
 \begin{aligned}
  \mU^{(1)} & = \nabla^a\left( \nabla^b(W_{acbd}W^{cefg}W^d{}_{efg}) + 2C_{cad}W^{cefg}W^d{}_{efg} \right) , \\
  \mU^{(2)} & = \nabla^a\Bigl( \nabla^b( W_{acde}W_b{}^c{}_{fg}W^{defg}) - 2W_{acde}W^{defg}C_{fg}{}^c \\
   & \qquad\qquad  - \frac{1}{6}\nabla_a( W_{bc}{}^{de}W_{de}{}^{fg}W_{fg}{}^{bc}) \Bigr) ,
 \end{aligned}
\end{equation}
where $C_{abc} := 2\nabla_{[a}P_{b]c}$ is the Cotton tensor.
Thus $\mU^{(1)}$ and $\mU^{(2)}$ are nontrivial scalar conformal invariants of weight $-8$ which are natural divergences in dimension $n=8$.

The remainder of this section is devoted to showing that any scalar conformal invariant of weight $-n$ on $n$-manifolds obtained by iteratively applying Proposition~\ref{prop:preserve-divergence} to a \ambientstraight\ symmetric tensor Riemannian invariant on ambient $(n+2)$-manifolds is a natural divergence.
In fact, we will prove a stronger result; see Proposition~\ref{prop:conf-inv-div}.
To that end, we first observe that the specific linear combination in Proposition~\ref{prop:preserve-divergence} is determined by requiring only that if $\cT_{A_1 \dotsm A_k}$ is annihilated by the infinitesimal generator $\cX := t\partial_t$ of dilation, then so too is $\cU_{A_1 \dotsm A_{k-1}}$.

\begin{lemma}
 \label{lem:ambient-recursive-construction}
 Let $\cT_{A_1 \dotsm A_k}$ be a symmetric tensor of weight $w$ on ambient $(n+2)$-manifolds such that $\cX^B\cT_{A_1 \dotsm A_{k-1}B} = 0$.
 Then the symmetric tensor
 \begin{equation*}
  \cU_{A_1 \dotsm A_{k-1}} := \cnabla^B \cT_{A_1 \dotsm A_{k-1} B} + \frac{k-1}{w-2k+2}\cnabla_{(A_1}\cT_{A_2 \dotsm A_{k-1})B}{}^B
 \end{equation*}
 of weight $w-2$ on ambient $(n+2)$-manifolds is such that $\cX^B \cU_{A_1 \dotsm A_{k-2} B} = 0$.
\end{lemma}

\begin{proof}
 Recall~\cite{GJMS1992}*{Equation~(1.4)} that $\cnabla\cX=\cg$.
 Since $\cT$ has rank $k$ and weight $w$, we compute that
 \begin{equation*}
  w\cT = L_{\cX}\cT = \cnabla_{\cX}\cT + k\cT .
 \end{equation*}
 Therefore $\cnabla_{\cX}\cT = (w-k)\cT$.
 Direct computation then gives
 \begin{align*}
  \cX^B \cU_{A_1 \dotsm A_{k-2} B} & = \cX^C \cnabla^B \cT_{A_1 \dotsm A_{k-2} CB} + \frac{1}{w-2k+2}\cX^C\cnabla_{C}\cT_{A_1 \dotsm A_{k-2}B}{}^B \\
   & \quad + \frac{k-2}{w-2k+2}\cX^C\cnabla_{(A_1}\cT_{A_2 \dotsm A_{k-2})CB}{}^B \\
  & = \left( -1 + \frac{w-k}{w-2k+2} - \frac{k-2}{w-2k+2} \right) \cT_{A_1 \dotsm A_{k-2} B}{}^B \\
  & = 0 . \qedhere 
 \end{align*}
\end{proof}

A straightforward computation implies that if $\cT_A$ is a tensor Riemannian invariant of weight $2-n$ on ambient $(n+2)$-manifolds that is annihilated by $\cX^A$, then $i^\ast\cnabla^A\cT_A$ is a scalar conformal invariant of weight $-n$ on $n$-manifolds that is a natural divergence, generalizing the invariants in Display~\eqref{eqn:dim8-examples}.
More generally:

\begin{proposition}
 \label{prop:conf-inv-div}
 Let $n \in \bN$ be even.
 Let $\cT_A$ be a tensor Riemannian invariant of rank $1$ and weight $-2k \geq 2-n$ on ambient $(n+2)$-manifolds.
 If $\cX^A\cT_A=0$, then the conformal invariant $\mI := i^\ast \cDelta^{n/2-k-1}\cnabla^A\cT_A$ of weight $-n$ on $n$-manifolds is a natural divergence.
\end{proposition}

\begin{remark}
 The following proof is similar in spirit to Marugame's proof~\cite{Marugame2017} that the $Q$-curvature of a pseudohermitian manifold is a divergence.
 Our proof follows Case and Yan's proof~\cite{CaseYan2024} of the formal self-adjointness of a family of conformally covariant differential operators, which include the GJMS operators~\cite{GJMS1992}.
\end{remark}

\begin{proof}
 Let $(\cmG,\cg)$ be a straight and normal ambient space for $(M,g)$.
 Set $r := \sqrt{-2\rho}$ and $s := rt$ in the domain $\cmG_+ := \{ \rho<0 \} \subset \cmG$.
 In these coordinates
 \begin{equation}
  \label{eqn:ambient-to-poincare-coordinates}
  \begin{aligned}
  \cg & = -ds^2 + s^2g_+ , \\
  g_+ & = r^{-2}(dr^2 + g_{-r^2/2}) ,
  \end{aligned}
 \end{equation}
 and $\cX_A = -s\,ds$.
 Consider the hypersurface $\hmG := \{ s=1 \} \subseteq \cmG_+$.
 On the one hand, since $\cT_A$ has weight $-2k$ and is annihilated by $\cX^A$, we see that $\cT = s^{-2k}\hpi^\ast\hT$ for some one-form $\hT$ on $\hmG$, where $\hpi \colon \cmG \to \hmG$ is the projection $\hpi(s,x,r) := (x,r)$.
 Therefore
 \begin{equation*}
  \cnabla^A\cT_A = s^{-2k-2}\hpi^\ast\hnabla^a\hT_a ,
 \end{equation*}
 where $\hnabla$ is the Levi-Civita connection of $g_+$.
 On the other hand, Equation~\eqref{eqn:ambient-to-poincare-coordinates} implies that
 \begin{equation*}
  \cDelta(s^w\hpi^\ast\hu) = s^{w-2}\hpi^\ast\bigl( \hDelta - w(n+w) \bigr)\hu
 \end{equation*}
 for all $w \in \bR$ and $\hu \in C^\infty(\hmG)$, where $\hDelta$ is the Laplacian of $g_+$.
 Therefore
 \begin{equation}
  \label{eqn:poincare-formula}
  \begin{aligned}
   \cI & := \cDelta^{n/2-k-1}\cnabla^A\cT_A \\
    & = s^{-n}\hpi^\ast\bigl(\hDelta + 2(n-2)\bigr)\dotsm\bigl(\hDelta + 2(k+1)(n-2k-2)\bigr)\hnabla^a\hT_a .
  \end{aligned}
 \end{equation}
 
 Consider now the integral
 \begin{equation*}
  \int_{r>\varepsilon} \cI \rv_{\hmG} \dvol_{g_+} 
 \end{equation*}
 for $\varepsilon>0$.
 On the one hand, the fact that $\cI \in \cmE[-n]$ and $i^\ast\cI = \mI$ implies that
 \begin{equation*}
  \cI\rv_{\hmG} = \mI r^n + O(r^{n+2})
 \end{equation*}
 as $r \to 0$.
 Since $g_+ = r^{-2}(dr^2 + g) + O(1)$, we deduce that
 \begin{equation}
  \label{eqn:log-part}
  \int_{r>\varepsilon} \cI\rv_{\hmG} \dvol_{g_+} = -\left(\int_M \mI^g \dvol_g \right)\log\varepsilon + O(1)
 \end{equation}
 as $\varepsilon\to0$.
 On the other hand, Equation~\eqref{eqn:poincare-formula} implies that
 \begin{equation*}
  \int_{r>\varepsilon} \cI\rv_{\hmG} \dvol_{g_+} = \int_{r>\varepsilon} \bigl( \hDelta + 2(n-2) \bigr) \dotsm \bigl( \hDelta + 2(k+1)(n-2k-2) \bigr)\hnabla^a\hT_a \dvol_{g_+} .
 \end{equation*}
 Since the right-hand side is a natural divergence, it can be written as the integral over $\{ r = \varepsilon \}$ of the normal component of an explicit natural vector field evaluated at $g_+ = \cg\rv_{T\hmG}$ times the volume element $\varepsilon^{-n}\dvol_{g_{-\varepsilon^2/2}}$.
 In particular, since $\cg$ is smooth, there can be no log-terms in the asymptotic expansion as $\varepsilon\to0$.
 Comparing this with Equation~\eqref{eqn:log-part} implies that $\int\mI\dvol=0$ on any compact Riemannian $n$-manifold.
 
 Finally, that $\int\mI\dvol=0$ for any metric on $T^n$ implies~\cite{Gilkey1975}*{Theorem~1.5(d)} that there is a constant $c \in \bR$ such that $\mI \equiv c\Pf$ modulo a natural divergence.
 Since also $\int\mI\dvol = 0$ on manifolds with nonzero Euler characteristic, which exist since $n$ is even, we conclude that $\mI$ is a natural divergence.
\end{proof}

\section{Proof of Theorem~\ref{thm:main-thm}}
\label{sec:cvi}

We begin by proving that the renormalized integral of a natural divergence vanishes.
We do so by adapting an argument used by Albin~\cite{Albin2009}*{Sections~3 and~4} to study the asymptotic expansions of scalar Riemannian invariants on Poincar\'e--Einstein manifolds.

\begin{lemma}
 \label{lem:renormalized-divergence-is-zero}
 Fix an even integer $n$.
 Let $F$ be a tensor Riemannian invariant of rank $1$ and weight at least $2-n$.
 If $(M^n,g_+)$ is a Poincar\'e--Einstein manifold, then
 \begin{equation*}
  \Rint \divsymb F \dvol = 0 .
 \end{equation*}
\end{lemma}

\begin{proof}
 Pick a representative $h$ of the conformal boundary of $(M^n,g_+)$.
 Chrus\'ciel, Delay, Lee, and Skinner showed~\cite{ChruscielDelayLeeSkinner2005}*{Theorem~A} that there is a collar diffeomorphism $\Phi$ from $[0,\varepsilon_0)_r \times \partial M$ onto a neighborhood of $\partial M$ in $M$ such that $\Phi(0,\cdot)$ is the inclusion $\partial M \hookrightarrow M$ and $\Phi^\ast g_+ = r^{-2}(dr^2 + g_r)$ for a one-parameter family $g_r$ of metrics on $\partial M$ such that
 \begin{equation}
  \label{eqn:pe-expansion}
  g_r = h + h_{(2)}r^2 + \dotsm + h_{(n-2)}r^{n-2} + kr^{n-1} + O(r^n) 
 \end{equation}
 has an even expansion to order $O(r^{n-1})$ and $k$ is trace-free with respect to $h$.
 
 Let $\varepsilon \in ( 0 , \varepsilon_0 )$.
 The Divergence Theorem implies that
 \begin{equation}
  \label{eqn:divergence-theorem}
  \int_{r>\varepsilon} (\divsymb F)^{g_+} \dvol_{g_+} = \oint_{r=\varepsilon} \lp F , \onf \rp^{g_+} \darea_{g_+} ,
 \end{equation}
 where $\onf$ is the outward-pointing unit normal and $\darea_{g_+}$ is the induced volume form on the level set $\{ r=\varepsilon \}$.
 Note that
 \begin{equation*}
  \darea_{g_+} = \varepsilon^{1-n}\dvol_{g_{\varepsilon}} .
 \end{equation*}
 Since $\tr_{h}k=0$, we deduce from Equation~\eqref{eqn:pe-expansion} that the asymptotic expansion
 \begin{equation}
  \label{eqn:area-expansion}
  \darea_{g_+} = \left( a_0\varepsilon^{1-n} + a_2\varepsilon^{3-n} + \dotsm + a_{n-2}\varepsilon^{-1} + a_n\varepsilon + O(\varepsilon^2) \right) \darea_h
 \end{equation}
 is odd to order $O(\varepsilon^2)$.
 We will show that the asymptotic expansion
 \begin{equation}
  \label{eqn:Fn-expansion}
  \lp F, \onf \rp^{g_+} = f_0 + f_2\varepsilon^2 + \dotsm + f_n\varepsilon^n + O(\varepsilon^{n+1})
 \end{equation}
 is even to order $O(\varepsilon^{n+1})$.
 Combining Equations~\eqref{eqn:divergence-theorem}, \eqref{eqn:area-expansion}, and~\eqref{eqn:Fn-expansion} then gives the final conclusion.
 
 Pick a point $p \in \partial M$.
 Let $\{ x^\alpha \}_{\alpha=1}^{n-1}$ be geodesic normal coordinates for $h$ near $p$ in $\partial M$.
 Extend these to coordinates on a neighborhood of $p$ in $M$ via the collar diffeomorphism $\Phi$ and let $\{ \partial_r , \partial_{x^\alpha} \}$ denote the resulting coordinate vector fields.
 Set $Y_0 := r\partial_r$ and $Y_\alpha := r\partial_{x^\alpha}$ for $\alpha \in \{ 1, \dotsc, n-1 \}$.
 Unless otherwise specified, in the remainder of this proof we use $a \in \{ 0, \dotsc, n-1 \}$ and $\alpha \in \{ 1, \dotsc, n-1 \}$ to denote components of Riemannian tensors for $g_+$ defined with respect to the frame $\{ Y_0, Y_\alpha \}$.
 For example, $(g_+)_{00} = 1$, $(g_+)_{0\alpha} = 0$, and $(g_+)_{\alpha\beta} = g_{\alpha\beta}$, where $g_{\alpha\beta}$ are the components of $g := r^2g_+$ with respect to $\{ \partial_r, \partial_{x^\alpha} \}$.
 Thus $(g_+^{-1})^{00}=1$, $(g_+^{-1})^{0\alpha}=0$, and $(g_+^{-1})^{\alpha\beta} = g^{\alpha\beta}$, where $g^{\alpha\beta}$ are the components of $g^{-1}$ with respect to $\{ \partial_r, \partial_{x^\alpha} \}$.
 
 Define $\gamma_{ab}^c$ by $\nabla^{g_+}_{Y_a}Y_b = \gamma_{ab}^cY_c$.
 Albin showed~\cite{Albin2009}*{Equation~(3.5)} that
 \begin{equation*}
  \gamma_{ab}^c = r\Gamma_{ab}^c - \delta_b^0\delta_a^c + \delta_0^cg_{ab} ,
 \end{equation*}
 where $\Gamma_{ab}^c$ denotes the Christoffel symbols of $g$.
 In particular,
 \begin{equation}
  \label{eqn:albin-connection}
  \begin{aligned}
   \gamma_{\alpha\beta}^\epsilon & = r\Gamma_{\alpha\beta}^\epsilon , & \gamma_{\alpha\beta}^0 & = -\frac{1}{2}r\partial_rg_{\alpha\beta} + g_{\alpha\beta} , \\
   \gamma_{\alpha0}^\epsilon & = \frac{1}{2}rg^{\epsilon\zeta}\partial_rg_{\alpha\zeta} - \delta_\alpha^\epsilon , & \gamma_{0\alpha}^\epsilon & = \frac{1}{2}rg^{\epsilon\zeta}\partial_rg_{\alpha\zeta} ,
  \end{aligned}
 \end{equation}
 and all other components of $\gamma_{ab}^c$ vanish.
 Denote by $\sF((-1)^\ell)$ the set of all functions $u \in C^\infty(X)$ such that the asymptotic expansion
 \begin{equation*}
  r^{(-1 + (-1)^\ell)/2}u = u_{(0)} + u_{(2)}r^{2} + \dotsm + u_{(n-2)}r^{n-2} + u_{(n-1)}r^{n-1} + O(r^{n})
 \end{equation*}
 is even to order $O(r^{n-1})$.
 That is, $\sF(1)$ consists of the functions that have an even expansion in $r$ to order $O(r^{n-1})$, and $\sF(-1)$ consists of the functions that have an odd expansion in $r$ to order $O(r^n)$.
 Note that if $u \in \sF((-1)^\ell)$ and $v \in \sF((-1)^m)$, then $uv \in \sF((-1)^{\ell+m})$.
 Since $g_{ab},g^{ab} \in \sF(1)$, we may freely raise and lower indices within $\sF((-1)^\ell)$.
 Equations~\eqref{eqn:albin-connection} imply that $\gamma_{ab}^c \in \sF((-1)^{\#_T(a,b,c)})$, where
 \begin{equation*}
  \#_T e := \left| \left\{ 1 \leq j \leq k \suchthat e_j \not= 0 \right\} \right|
 \end{equation*}
 denotes the number of nonzero terms in the multi-index $e \in \{ 0, \dotsc, n-1 \}^k$.
 
 Consider now the components $R_{ab}{}^c{}_d$ of the Riemann curvature tensor of $g_+$,
 \begin{equation*}
  R_{ab}{}^c{}_dY_c = \nabla^{g_+}_{Y_a}\nabla^{g_+}_{Y_b}Y_d - \nabla^{g_+}_{Y_b}\nabla^{g_+}_{Y_a}Y_d - \nabla^{g_+}_{[Y_a,Y_b]}Y_d .
 \end{equation*}
 Albin showed~\cite{Albin2009}*{Equation~(3.6)} that
 \begin{equation}
  \label{eqn:albin-curvature}
  R_{ab}{}^c{}_d = \gamma_{ae}^{c}\gamma_{bd}^{e} - \gamma_{be}^{c}\gamma_{ad}^{e} + Y_a(\gamma_{bd}^{c}) - Y_b(\gamma_{ad}^{c}) - \delta_a^0\gamma_{bd}^{c} + \delta_b^0\gamma_{ad}^{c} .
 \end{equation}
 Since $Y_0\bigl(\sF((-1)^m)\bigr) \subseteq \sF((-1)^m)$ and $Y_\alpha\bigl(\sF((-1)^m)\bigr) \subseteq \sF((-1)^{m+1})$, we deduce that $R_{abcd} \in \sF((-1)^{\#_T(a,b,c,d)})$.
 
 Denote the components of $\nabla^\ell\Rm$ by
 \begin{multline}
  \label{eqn:curvature-derivatives}
  R_{abcd;e_1\dotsm e_\ell} := Y_{e_\ell}(R_{abcd;e_1\dotsm e_{\ell-1}}) - \gamma_{e_\ell a}^fR_{fbcd;e_1\dotsm e_{\ell-1}} \\- \dotsm - \gamma_{e_\ell e_{\ell-1}}^f R_{abcd;e_1 \dotsm e_{\ell-2}f} .
 \end{multline}
 It follows from induction that
 \begin{equation*}
  R_{abcd;e_1 \dotsm e_{\ell}} \in \sF\bigl( (-1)^{\#_T(a,b,c,d,e)} \bigr) .
 \end{equation*}
 
 Suppose first that $F$ is odd.
 Then the Hodge dual $\ast F$ is an even natural $(n-1)$-form of weight at least $0$.
 Since $n-1$ is odd, a result of Gilkey~\cite{Gilkey1975}*{Theorem~(1.2)} implies that $F=0$.
 Hence we may assume that $F$ is even.
 Then $\lp F, \onf \rp$ is a linear combination of complete contractions of
 \begin{equation}
  \label{eqn:product}
  \nabla^{N_1}\Rm \otimes \dotsm \otimes \nabla^{N_p}\Rm \otimes \mathop{\onf} ,
 \end{equation}
 where $N_1,\dotsc,N_p$ denote powers.
 Since $\onf = -Y_0$ along $\partial M$, any such complete contraction is a linear combination of products of $R_{abcd;e_1 \dotsm e_\ell}$, and each product has an even number of tangential indices.
 From the previous discussion, we conclude that $\lp F, \onf \rp^{g_+} \in \sF(1)$.
 Moreover, the coefficient of $r^{n-1}$ of $\lp F, \onf \rp^{g_+}$ is computed by multiplying the coefficient of $r^{n-1}$ in one factor $R_{abcd;e_1\dotsm e_\ell}$ by the coefficients of $r^0$ in the remaining factors.
 
 Given $v \in C^\infty(M)$, denote by $v^{(0)}$ and $v^{(n-1)}$ the coefficients of $r^0$ and $r^{n-1}$, respectively, in the Taylor expansion of $v$ near $\partial M$.
 Equations~\eqref{eqn:albin-connection} and~\eqref{eqn:albin-curvature} yield
 \begin{align*}
  \bigl(\gamma_{\alpha\beta}^0\bigr)^{(0)} & = h_{\alpha\beta} , & \bigl(\gamma_{\alpha\beta}^0\bigr)^{(n-1)} & = -\frac{n-3}{2}k_{\alpha\beta} , \\
  \bigl(\gamma_{\alpha0}^\beta\bigr)^{(0)} & = -\delta_\alpha^\beta , & \bigl(\gamma_{\alpha0}^\beta\bigr)^{(n-1)} & = \frac{n-1}{2}h^{\beta\epsilon}k_{\alpha\epsilon} , \\
  \bigl(\gamma_{0\alpha}^\beta\bigr)^{(0)} & = 0 , & \bigl(\gamma_{0\alpha}^\beta\bigr)^{(n-1)} & = \frac{n-1}{2}h^{\beta\epsilon}k_{\alpha\epsilon} , \\
  \bigl(R_{\alpha\beta\epsilon\zeta}\bigr)^{(0)} & = -2h_{\alpha[\epsilon}h_{\zeta]\beta} , & \bigl(R_{\alpha\beta\epsilon\zeta}\bigr)^{(n-1)} & = (n-3)(h_{\alpha[\epsilon}k_{\zeta]\beta} - h_{\beta[\epsilon}k_{\zeta]\alpha}) , \\
  \bigl(R_{\alpha0\beta0}\bigr)^{(0)} & = -h_{\alpha\beta} , & \bigl(R_{\alpha0\beta0}\bigr)^{(n-1)} & = -\frac{n^2-4n+5}{2}k_{\alpha\beta} .
 \end{align*}
 Combining this with Equation~\eqref{eqn:curvature-derivatives}, the previous paragraph, and naturality implies that the coefficient of $r^{n-1}$ of $\lp F,\onf\rp^{g_+}$ must be proportional to $\tr_{h}k = 0$.
\end{proof}

Combining Proposition~\ref{prop:easy-straight-derivatives} and Lemma~\ref{lem:renormalized-divergence-is-zero} yields our formula for the renormalized integral of a \straight\ scalar Riemannian invariant.

\begin{proof}[Proof of Theorem~\ref{thm:main-thm}]
 Let $I$ be a \straight\ scalar Riemannian invariant of weight $-2k \geq -n$.
 If $I$ is odd, then a parity argument~\citelist{ \cite{BaileyEastwoodGraham1994}*{p.\ 502} \cite{Gilkey1975}*{Theorem~1.2} } implies that $I$ is the Hodge dual of a top-degree Pontrjagin form.
 Such forms are conformally invariant~\cites{Avez1970}, and hence $\smallRint I \dvol = \int_M I\dvol$ in this case.
 We may thus assume that $I$ is even.

 Set $\ell := n/2-k$ and let $I_{\ell}$ be the \straight\ scalar Riemannian invariant determined by Proposition~\ref{prop:easy-straight-derivatives}.
 Let $\mI_\ell$ be the conformal invariant determined ambiently by $\cDelta^\ell\cI$, where $\cI$ is a \ambientstraight\ scalar Riemannian invariant to which $I$ is associated.
 On the one hand, since $\mI_\ell$ has weight $-n$, we see that
 \begin{equation*}
  \Rint \mI_\ell \dvol = \int_M \mI_\ell \dvol .
 \end{equation*}
 On the other hand, Lemma~\ref{lem:straight-is-almost-conformal}, Equation~\eqref{eqn:mI-constant-term}, and Lemma~\ref{lem:renormalized-divergence-is-zero} imply that
 \begin{equation*}
  \Rint \mI_\ell \dvol = \Rint I_\ell \dvol = 2^{n/2-k}\frac{(n/2-1)!(n-2k-1)!!}{(k-1)!}\Rint I \dvol . \qedhere
 \end{equation*}
\end{proof}

\section{Applications}
\label{sec:algorithm}

Theorem~\ref{thm:main-thm} and Lemma~\ref{lem:weyl-is-straight} together compute the renormalized integral of any scalar polynomial in the Weyl tensor.
This allows us to derive an explicit formula expressing the renormalized volume of an even-dimensional Poincar\'e--Einstein manifold in terms of its Euler characteristic and the convergent integral of a scalar conformal invariant.
The key observation is that the Pfaffian $\Pf = \Pf_{n/2}(\Rm)$ of an Einstein manifold is a scalar polynomial in the Weyl tensor (cf.\ \cite{Albin2009}*{Lemma~4.4}).

\begin{lemma}
 \label{lem:einstein-pfaffian}
 Let $(M^n,g)$ be an even-dimensional Einstein manifold.
 Then
 \begin{equation*}
  \Pf = \sum_{\ell=0}^{n/2} (n-2\ell-1)!! \left(\frac{2}{n}\trSch\right)^{n/2-\ell} \Pf_\ell(\Weyl) .
 \end{equation*}
\end{lemma}

\begin{proof}
 The assumption that $(M^n,g)$ is Einstein implies that
 \begin{equation*}
  R_{a_1a_2}{}^{b_1b_2} = \Weyl_{a_1a_2}{}^{b_1b_2} + \frac{4}{n}\trSch\delta_{a_1a_2}^{b_1b_2} .
 \end{equation*}
 We deduce from Equation~\eqref{eqn:generalized-pfaffian} that
 \begin{multline*}
  \Pf = 2^{-n/2}(n-1)!!\sum_{\ell=0}^{n/2} \binom{n/2}{\ell} \left(\frac{4}{n}\trSch\right)^{n/2-\ell} \\
   \times \delta_{b_1 \dotsm b_n}^{a_1 \dotsm a_n}\Weyl_{a_1a_2}{}^{b_1b_2} \dotsm \Weyl_{a_{2\ell-1}a_{2\ell}}{}^{b_{2\ell-1}b_{2\ell}} \delta_{a_{2\ell+1}}^{b_{2\ell+1}} \dotsm \delta_{a_n}^{b_n} .
 \end{multline*}
 The conclusion follows from the identity
 \begin{equation*}
  \delta_{b_1 \dotsm b_k}^{a_1 \dotsm a_k} \delta_{a_k}^{b_k} = \frac{n-k+1}{k}\delta_{b_1 \dotsm b_{k-1}}^{a_1 \dotsm a_{k-1}} . \qedhere
 \end{equation*}
\end{proof}

Since $\Pf_{\ell}(\Weyl)$ is a \straight\ scalar Riemannian invariant of weight $-2\ell$, we can apply Theorem~\ref{thm:main-thm} to produce an explicit ambient formula for the scalar conformal invariant in Equation~\eqref{eqn:renormalized-volume-cqy}.
More generally, we produce a local formula for the Pfaffian of an Einstein manifold modulo natural divergences.

\begin{corollary}
 \label{cor:pfaffian-via-straight}
 Let $(M^n,g)$ be an even-dimensional Einstein manifold.
 Then
 \begin{equation*}
  \Pf \equiv (n-1)!!\left( \frac{2}{n}\trSch \right)^{n/2} + \sum_{\ell=2}^{n/2} \left(-2\right)^{\ell-n/2}\frac{(\ell-1)!}{(n/2-1)!}\mP_{\ell,n} \mod \im \Delta ,
 \end{equation*}
 where $\mP_{\ell,n} := i^\ast\left( \cDelta^{n/2-\ell}\Pf_\ell\bigl(\cRm\bigr) \right)$.
\end{corollary}

\begin{proof}
  On the one hand, Lemma~\ref{lem:einstein-pfaffian} and the fact $\Pf_1(W)=0$ yield
  \begin{equation*}
    \Pf = (n-1)!!\left( \frac{2}{n} \trSch \right)^{n/2} + \sum_{\ell=2}^{n/2} (n-2\ell-1)!!\left( \frac{2}{n} \trSch \right)^{n/2-\ell}\Pf_\ell(\Weyl) .
  \end{equation*}
  On the other hand, Proposition~\ref{prop:easy-straight-derivatives} yields
  \begin{equation*}
    (\Pf_{\ell}(\Weyl))_{n/2-\ell} \equiv \left( -\frac{4}{n} \trSch \right)^{n/2-\ell} \frac{(n/2-1)!(n-2\ell-1)!!}{(\ell-1)!}\Pf_\ell(\Weyl) \mod \im\Delta .
  \end{equation*}
  Combining these with Lemma~\ref{lem:straight-is-almost-conformal} yields the conclusion.
\end{proof}

Corollary~\ref{cor:pfaffian-via-straight} yields our explicit formula for the (renormalized) volume of an even-dimensional (Poincar\'e--)Einstein manifold.

\begin{proof}[Proof of Theorem~\ref{thm:cqy}]
 Apply Equation~\eqref{eqn:renormalized-pfaffian-formula}, Lemma~\ref{lem:renormalized-divergence-is-zero}, and Corollary~\ref{cor:pfaffian-via-straight}.
\end{proof}

\begin{proof}[Proof of Theorem~\ref{thm:gbc}]
 Apply Corollary~\ref{cor:pfaffian-via-straight}.
\end{proof}

By direct computation, any complete contraction of $\Weyl^{\otimes2}$ is proportional to
\begin{subequations}
 \label{eqn:weyl-contractions}
 \begin{equation}
  \label{eqn:weyl-contractions2}
  \mW_{2,1} := \Weyl_{abcd}\Weyl^{abcd} ;
 \end{equation}
 any complete contraction of $\Weyl^{\otimes 3}$ is a linear combination of
 \begin{equation}
  \label{eqn:weyl-contractions3}
  \begin{aligned}
   \mW_{3,1} & := \Weyl_{ab}{}^{cd}\Weyl_{cd}{}^{ef}\Weyl_{ef}{}^{ab} , & \mW_{3,2} & := \Weyl_a{}^c{}_b{}^d \Weyl_c{}^e{}_d{}^f \Weyl_e{}^a{}_f{}^b ;
  \end{aligned}
 \end{equation}
 and any complete contraction of $\Weyl^{\otimes 4}$ is a linear combination of
 \begin{equation}
  \label{eqn:weyl-contractions4}
  \begin{aligned}
   \mW_{4,1} & := \Weyl_{abcd}\Weyl^{abcd}\Weyl_{efgh}\Weyl^{efgh} , & \mW_{4,2} & := \Weyl_{ab}{}^{cd}\Weyl_{cd}{}^{ef}\Weyl_{ef}{}^{gh}\Weyl_{gh}{}^{ab} , \\
   \mW_{4,3} & := \Weyl_{acde}\Weyl_b{}^{cde}\Weyl^a{}_{fgh}\Weyl^{bfgh} , & \mW_{4,4} & := \Weyl_{abcd}\Weyl^{cd}{}_{ef}\Weyl^a{}_g{}^e{}_h\Weyl^{bgfh} , \\
   \mW_{4,5} & := \Weyl_{abcd}\Weyl^{cd}{}_{ef}\Weyl^{ae}{}_{gh}\Weyl^{bfgh} , & \mW_{4,6} & := \Weyl_a{}^c{}_b{}^d\Weyl_c{}^e{}_d{}^f\Weyl_e{}^g{}_f{}^h\Weyl_g{}^a{}_h{}^b , \\
   \mW_{4,7} & := \Weyl_a{}^c{}_b{}^d\Weyl_{ecfd}\Weyl^a{}_g{}^e{}_h \Weyl^{bgfh} ;
  \end{aligned}
 \end{equation}
\end{subequations}
see~\cite{FullingKingWybourneCummins1992} for details.
We express $\Pf_{\ell}(\Weyl)$, $\ell \leq 4$, in terms of these contractions:

\begin{lemma}
 \label{lem:low-order-weyl-pfaffian}
 Let $(M^n,g)$ be a Riemannian manifold.
 Then
 \begin{align*}
  \Pf_2(\Weyl) & = \frac{1}{8}\mW_{2,1} , \\
  \Pf_3(\Weyl) & = \frac{1}{12}\mW_{3,1} - \frac{1}{6}\mW_{3,2} , \\
  \Pf_4(\Weyl) & = \frac{1}{128}\mW_{4,1} + \frac{1}{64}\mW_{4,2} - \frac{1}{8}\mW_{4,3} - \frac{1}{4}\mW_{4,4} + \frac{1}{8}\mW_{4,5} + \frac{1}{8}\mW_{4,6} - \frac{1}{4}\mW_{4,7} .
 \end{align*}
\end{lemma}

\begin{proof}
 Recall that the Weyl tensor is totally trace-free.
 Direct computation from the definition of $\Pf_\ell$ yields
 \begin{align*}
  \Pf_2(\Weyl) & = \frac{1}{2^22!}\Weyl_{abcd}\Weyl^{abcd} , \\
  \Pf_3(\Weyl) & = \frac{1}{2^33!}\left( 2\Weyl_{ab}{}^{ef}\Weyl_{cd}{}^{ab}\Weyl_{ef}{}^{cd} - 8\Weyl_{ab}{}^{ce}\Weyl_{cd}{}^{af}\Weyl_{ef}{}^{bd} \right) , \\
  \Pf_4(\Weyl) & = \frac{1}{2^44!}\Bigl( 3\Weyl_{ab}{}^{cd}\Weyl_{cd}{}^{ab}\Weyl_{ef}{}^{gh}\Weyl_{gh}{}^{ef} + 6\Weyl_{ab}{}^{cd}\Weyl_{cd}{}^{ef}\Weyl_{ef}{}^{gh}\Weyl_{gh}{}^{ab} \\
   & \qquad - 48\Weyl_{ab}{}^{cg}\Weyl_{cd}{}^{ab}\Weyl_{ef}{}^{dh}\Weyl_{gh}{}^{ef} - 96\Weyl_{ab}{}^{eg}\Weyl_{cd}{}^{ab}\Weyl_{ef}{}^{ch}\Weyl_{gh}{}^{df} \\
   & \qquad + 48\Weyl_{ab}{}^{ce}\Weyl_{cd}{}^{ag}\Weyl_{ef}{}^{bh}\Weyl_{gh}{}^{df} - 96\Weyl_{ab}{}^{eg}\Weyl_{cd}{}^{ah}\Weyl_{ef}{}^{bc}\Weyl_{gh}{}^{df} \Bigr) .
 \end{align*}
 The conclusion now follows from the Bianchi identity $\Weyl_{[abc]d}=0$:
 \begin{align*}
  \Weyl_{ab}{}^{ce}\Weyl_{cd}{}^{af}\Weyl_{ef}{}^{bd} & = \Weyl_{abc}{}^e\Weyl^{afcd}\Weyl^b{}_{fed} - \Weyl_{abc}{}^e\Weyl^{afcd}\Weyl^b{}_{efd} \\
   & = \Weyl_{abc}{}^e\Weyl^{afcd}\Weyl^b{}_{fed} - \frac{1}{4}\Weyl_{acb}{}^e\Weyl^{acfd}\Weyl^b{}_{efd}
 \end{align*}
 and
 \begin{align*}
  \MoveEqLeft \Weyl_{ab}{}^{eg}\Weyl_{cd}{}^{ab}\Weyl_{ef}{}^{ch}\Weyl_{gh}{}^{df} = \Weyl_{cd}{}^{ab}\Weyl_{abeg}\Weyl^{chef}\Weyl^d{}_f{}^g{}_h \\
   & = \Weyl_{cd}{}^{ab}\Weyl_{abeg}\Weyl^{chef}\Weyl^d{}_h{}^g{}_f + \Weyl_{cd}{}^{ab}\Weyl_{abeg}\Weyl^{chef}\Weyl^{dg}{}_{fh} \\
   & = \Weyl_{cd}{}^{ab}\Weyl_{abeg}\Weyl^{chef}\Weyl^d{}_h{}^g{}_f - \frac{1}{2}\Weyl_{cd}{}^{ab}\Weyl_{abeg}\Weyl^{cefh}\Weyl^{dg}{}_{fh} . \qedhere
 \end{align*}
\end{proof}

Lemma~\ref{lem:low-order-weyl-pfaffian} allows us to compute $\mP_{\ell,n}$ for all integers $2 \leq \ell \leq 4$.

\begin{corollary}
 \label{cor:mPln-low-dimensions}
 Let $(M^n,g)$ be an even-dimensional Riemannian manifold.  Then
 \begin{align*}
  \mP_{2,n} & = \frac{1}{8}i^\ast\bigl( \cDelta^{n/2-2}\cmW_{2,1} \bigr) , && \text{if $n \geq 4$},  \\
  \mP_{3,n} & = \frac{1}{12}i^\ast\bigl( \cDelta^{n/2-3}( \cmW_{3,1} - 2\cmW_{3,2} ) \bigr) , && \text{if $n \geq 6$}, \\
  \mP_{4,n} & = \frac{1}{128}i^\ast\bigl( \cDelta^{n/2-4}( \cmW_{4,1} + 2\cmW_{4,2} - 16\cmW_{4,3} - 32\cmW_{4,4} \\
   & \qquad\qquad\quad+ 16\cmW_{4,5} + 16\cmW_{4,6} - 32\cmW_{4,7} ) \bigr) , && \text{if $n \geq 8$} ,
 \end{align*}
 where $\cmW_{j,k}$ are defined by the same contractions of $\cRm^{\otimes j}$ as in Equation~\eqref{eqn:weyl-contractions}.
\end{corollary}

\begin{proof}
 Combine Lemma~\ref{lem:low-order-weyl-pfaffian} with the definition of $\mP_{\ell,n}$ in Theorem~\ref{thm:cqy}.
\end{proof}

Corollary~\ref{cor:mPln-low-dimensions} allows us to express the Gauss--Bonnet formula of Theorem~\ref{thm:cqy} in a more familiar way in low dimensions.
For example, we immediately get the following formula in dimension eight:

\begin{corollary}
 \label{cor:gbc8}
 Let $(M^8,g_+)$ be an eight-dimensional Poincar\'e--Einstein manifold.
 Then
 \begin{multline*}
  16\pi^4\chi(M) = 105\mV + \frac{1}{192}\int_M \iota^\ast \biggl( \cDelta^2\cmW_{2,1} - \frac{8}{3}\cDelta\cmW_{3,1} + \frac{16}{3}\cDelta\cmW_{3,2} + \frac{3}{2}\cmW_{4,1} \\
   + 3\cmW_{4,2} - 24\cmW_{4,3} - 48\cmW_{4,4} + 24\cmW_{4,5} + 24\cmW_{4,6} - 48\cmW_{4,7} \biggr) \dvol .
 \end{multline*}
\end{corollary}

\begin{proof}
 Apply Theorem~\ref{thm:cqy} and Corollary~\ref{cor:mPln-low-dimensions}.
\end{proof}

\begin{remark}
 Theorem~\ref{thm:gbc} and Corollary~\ref{cor:mPln-low-dimensions} imply that if $(M^8,g)$ is a compact Einstein eight-manifold with $\Ric_g = 7\lambda g$, then 
 \begin{multline*}
  16\pi^4\chi(M) = 105\lambda^4\Vol_g(M) + \frac{1}{192}\int_M \iota^\ast \biggl( \cDelta^2\cmW_{2,1} - \frac{8}{3}\cDelta\cmW_{3,1} + \frac{16}{3}\cDelta\cmW_{3,2} + \frac{3}{2}\cmW_{4,1} \\
   + 3\cmW_{4,2} - 24\cmW_{4,3} - 48\cmW_{4,4} + 24\cmW_{4,5} + 24\cmW_{4,6} - 48\cmW_{4,7} \biggr) \dvol .
 \end{multline*}
\end{remark}

With a little more work, we also recover the Gauss--Bonnet formula of Chang, Qing, and Yang~\cite{ChangQingYang2006}*{Equation~(4.8)} in dimension six.

\begin{corollary}
 \label{cor:gbc6}
 Let $(M^6,g_+)$ be a six-dimensional Poincar\'e--Einstein manifold.
 Then
 \begin{equation*}
  8\pi^3\chi(M) = -15\mV - \frac{1}{16}\int_M \iota^\ast \left( \lv\cnabla\cRm\rv^2 - \frac{7}{3}\cmW_{3,1} - \frac{4}{3}\cmW_{3,2} \right) \dvol .
 \end{equation*}
\end{corollary}

\begin{proof}
 On the one hand, Theorem~\ref{thm:cqy} and Corollary~\ref{cor:mPln-low-dimensions} yield
 \begin{equation*}
  8\pi^3\chi(M) = -15\mV + \int_M \iota^\ast \left( -\frac{1}{32}\cDelta\cmW_{2,1} + \frac{1}{12}\cmW_{3,1} - \frac{1}{6}\cmW_{3,2} \right) \dvol .
 \end{equation*}
 On the other hand, the second Bianchi identity and the fact $\cRic=0$ imply that
 \begin{equation*}
  \cDelta\cmW_{2,1} = 2\lv\cnabla\cRm\rv^2 - 2\cmW_{3,1} - 8\cmW_{3,2} .
 \end{equation*}
 The conclusion readily follows.
\end{proof}

\begin{remark}
 A similar argument implies that if $(M^6,g)$ is a compact Einstein six-manifold with $\Ric_g = 5\lambda g$, then
 \begin{equation*}
  8\pi^3\chi(M) = 15\lambda^3\Vol_g(M) - \frac{1}{16}\int_M \iota^\ast \left( \lv\cnabla\cRm\rv^2 - \frac{7}{3}\cmW_{3,1} - \frac{4}{3}\cmW_{3,2} \right) \dvol .
 \end{equation*}
\end{remark}

We conclude this paper by computing three additional examples of renormalized curvature integrals which indicate the applicability of our results.

Our first example, when combined with the second Bianchi identity, computes the renormalized integral of any complete contraction of $(\nabla\Weyl)^{\otimes 2}$.

\begin{example}
 \label{ex:nablaW2}
 Let $(M^n,g_+)$, $n \geq 6$, be an even-dimensional Poincar\'e--Einstein manifold.
 Then
 \begin{equation*}
  \Rint \lv\nabla\Weyl\rv^2 \dvol = \frac{2^{3-n/2}2!}{(n/2-1)!(n-7)!!} \int_M i^\ast(\cDelta^{n/2-3}\cL) \dvol ,
 \end{equation*}
 where
 \begin{equation*}
  \cL := \frac{n-1}{2(n-5)}\cDelta\cmW_{2,1} + \cmW_{3,1} + 4\cmW_{3,2} .
 \end{equation*}
 To see this, recall that the Bianchi identity implies that if $\Ric_g = 2(n-1)\lambda g$, then
 \begin{equation}
  \label{eqn:DeltaWeyl}
  \Delta\Weyl_{abcd} = 4(n-1)\lambda\Weyl_{abcd} - \Weyl_{ab}{}^{ef}\Weyl_{efcd} - 2\Weyl_{aecf}\Weyl_b{}^e{}_d{}^f + 2\Weyl_{aedf}\Weyl_b{}^e{}_c{}^f .
 \end{equation}
 Now, Lemma~\ref{lem:renormalized-divergence-is-zero} yields
 \begin{equation*}
  \Rint \lv\nabla\Weyl\rv^2 \dvol = -\Rint \lp \Weyl, \Delta\Weyl \rp \dvol ,
 \end{equation*}
 while Equation~\eqref{eqn:DeltaWeyl} and the assumption $\Ric_{g_+} = -(n-1)g_+$ yield
 \begin{equation*}
  \Rint \lp \Weyl, \Delta\Weyl \rp \dvol = -\Rint \Bigl( 2(n-1)\mW_{2,1} + \mW_{3,1} + 4\mW_{3,2} \Bigr) \dvol .
 \end{equation*}
 The claim now follows from Theorem~\ref{thm:main-thm}.
 \qed
\end{example}

The explanation of our second example indicates how to compute the renormalized integral of a complete contraction of $(\nabla^2\Weyl)^{\otimes 2}$ or $\Weyl \otimes (\nabla\Weyl)^{\otimes 2}$.

\begin{example}
 Let $(M^n,g_+)$, $n \geq 8$, be an even-dimensional Poincar\'e--Einstein manifold.
 Then
 \begin{equation*}
  \Rint \lv\nabla^2\Weyl\rv^2 \dvol = \frac{2^{4-n/2}3!}{(n/2-1)!(n-9)!!} \int_M i^\ast(\cDelta^{n/2-4}\cL) \dvol  ,
 \end{equation*}
 where
 \begin{multline*}
  \cL := \frac{(n-1)(3n+1)}{12(n-5)(n-7)}\cDelta^2\cmW_{2,1} + \frac{5n+3}{6(n-7)}\cDelta\cmW_{3,1} + \frac{2(5n+3)}{3(n-7)}\cDelta\cmW_{3,2} \\
   + \cmW_{4,2} + 4\cmW_{4,3} + 8\cmW_{4,4} - 6\cmW_{4,5} + 8\cmW_{4,6} - 8\cmW_{4,7} .
 \end{multline*}
 Indeed, the Ricci identity and the assumption $\Ric_{g_+} = -(n-1)g_+$ yield
 \begin{multline*}
  \Delta\nabla_e\Weyl_{abcd} = \nabla_e\Delta\Weyl_{abcd} - (n+3)\nabla_e\Weyl_{abcd} \\
   - 4\Weyl_{feh}{}_{[a}\nabla^f\Weyl^h{}_{b]cd} - 4\Weyl_{feh}{}_{[c}\nabla^f\Weyl^h{}_{d]ab} .
 \end{multline*}
 Combining this with Lemma~\ref{lem:renormalized-divergence-is-zero} yields
 \begin{multline*}
  \Rint \lv\nabla^2\Weyl\rv^2 \dvol = \Rint \Bigl( \lv\Delta\Weyl\rv^2 - (n+3)\lp\Weyl,\Delta\Weyl\rp \\
   + 8\Weyl^{abcd}\Weyl_{efha}\nabla^e\nabla^f\Weyl^h{}_{bcd} \Bigr) \dvol .
 \end{multline*}
 On the one hand, the first Bianchi identity yields
 \begin{align*}
  \MoveEqLeft \Weyl_{aecf}\Weyl_b{}^e{}_d{}^f\Weyl^{agdh}\Weyl^b{}_g{}^c{}_h  = \Weyl_{aecf}\Weyl_b{}^e{}_d{}^f\Weyl^{agdh}(\Weyl^c{}_g{}^b{}_h + \Weyl^{bc}{}_{gh}) \\
   & = \Weyl_{aecf}(\Weyl_d{}^e{}_b{}^f + \Weyl_{bd}{}^{ef})\Weyl^{agdh}\Weyl^c{}_g{}^b{}_h + \frac{1}{2}\Weyl_{aecf}\Weyl_b{}^e{}_d{}^f\Weyl^{adgh}\Weyl^{bc}{}_{gh} \\
   & = \mW_{4,7} - \mW_{4,4} + \frac{1}{4}\mW_{4,5} .
 \end{align*}
 Combining this with Equation~\eqref{eqn:DeltaWeyl} yields
 \begin{align*}
  \MoveEqLeft[10] \Rint \left( \lv\Delta \Weyl \rv^2 - (n+3)\lp\Weyl,\Delta\Weyl\rp \right) \dvol = \Rint \Bigl( 2(n-1)(3n+1)\mW_{2,1} \\
   & + (5n-1)\mW_{3,1} + 4(5n-1)\mW_{3,2} \\
   & + \mW_{4,2} + 16\mW_{4,4} - 2\mW_{4,5} + 8\mW_{4,6} - 8\mW_{4,7} \Bigr) \dvol .
 \end{align*}
 On the other hand, the Ricci identity yields
 \begin{multline*}
  \Rint 8\Weyl^{abcd}\Weyl_{efha}\nabla^e\nabla^f\Weyl^h{}_{bcd} \dvol \\
   = \Rint 4\Bigl( \mW_{3,1} + 4\mW_{3,2} + \mW_{4,3} - 2\mW_{4,4} - \mW_{4,5} \Bigr) \dvol .
 \end{multline*}
 Combining these formulas with Theorem~\ref{thm:main-thm} yields the final claim.
 \qed
\end{example}

Our final example illustrates a more general approach to applying Propositions~\ref{prop:easy-straight-derivatives} and~\ref{prop:preserve-divergence} to compute renormalized curvature integrals.

\begin{example}
 Let $(M^n,g_+)$, $n \geq 10$, be an even-dimensional Poincar\'e--Einstein manifold.
 Then
 \begin{equation*}
  \Rint \left| \nabla\lv\Weyl\rv^2 \right|^2 \dvol = \frac{2^{5-n/2}4!}{(n/2-1)!(n-11)!!}\int_M i^\ast(\cDelta^{n/2-5}\cL) \dvol ,
 \end{equation*}
 where
 \begin{equation*}
  \cL := -\cmW_{2,1}\cDelta\cmW_{2,1} + \frac{(n-5)}{2(n-9)}\cDelta\cmW_{2,1}^2 .
 \end{equation*}
 Indeed, Lemma~\ref{lem:renormalized-divergence-is-zero} yields
 \begin{align*}
  \Rint \left| \nabla\lv\Weyl\rv^2 \right|^2 \dvol & = \Rint (-\lv\Weyl\rv^2\Delta\lv\Weyl\rv^2) \dvol \\
  & = \Rint \left( -\lv\Weyl\rv^2(\lv\Weyl\rv^2)_1 + 4(n-5)\lv\Weyl\rv^4 \right) \dvol ,
 \end{align*}
 where $(\lv\Weyl\rv^2)_1 := \Delta\lv\Weyl\rv^2 + 4(n-5)\lv\Weyl\rv^2$ is the specialization to $(M^n,g_+)$ of the \straight\ scalar Riemannian invariant produced by Proposition~\ref{prop:easy-straight-derivatives}.
 The conclusion now follows from Theorem~\ref{thm:main-thm}.
 \qed
\end{example}

The general principle is that Propositions~\ref{prop:easy-straight-derivatives} and~\ref{prop:preserve-divergence} allow one to write certain divergences as a linear combination of a \straight\ and a lower-order tensor Riemannian invariant.
We speculate that, when combined with the Bianchi identity, this procedure allows one to compute \emph{any} renormalized curvature integral.

\section*{Acknowledgements}
This project is a part of the AIM SQuaRE ``Global invariants of Poincar\'e--Einstein manifolds and applications''.
We thank the American Institute for Mathematics for their support.
We also thank Samuel Blitz, Danilo D\'iaz, Matthew Gursky, and the anonymous referee for valuable comments.

\section*{Funding}
JSC acknowledges support from a Simons Foundation Collaboration Grant for Mathematicians, ID 524601.
YJL acknowledges support from an NSF-LEAPS grant (Grant \#DMS-2418740).
WY acknowledges support from the NSFC (Grant No.\ 12071489, No.\ 12025109) and from the GuangDong Basic and Applied Basic Research Foundation (Grant No.\ 2023A1515012877).

\section*{Declarations}
The authors have no relevant financial or non-financial interests to disclose.

\section*{Data availability}
Data sharing not applicable to this article as no datasets were generated or analyzed during the study.

\bibliography{../Bibliography/bib}

\newcommand{\noopsort}[1]{}
\begin{bibdiv}
\begin{biblist}

\bib{Albin2009}{article}{
      author={Albin, Pierre},
       title={Renormalizing curvature integrals on {Poincar\'e}-{Einstein} manifolds},
        date={2009},
        ISSN={0001-8708},
     journal={Adv. Math.},
      volume={221},
      number={1},
       pages={140\ndash 169},
         url={http://dx.doi.org/10.1016/j.aim.2008.12.002},
      review={\MR{2509323}},
}

\bib{Alexakis2012}{book}{
      author={Alexakis, Spyros},
       title={The decomposition of global conformal invariants},
      series={Annals of Mathematics Studies},
   publisher={Princeton University Press},
     address={Princeton, NJ},
        date={2012},
      volume={182},
        ISBN={978-0-691-15348-3},
      review={\MR{2918125}},
}

\bib{Anderson2001b}{article}{
      author={Anderson, Michael~T.},
       title={{$L^2$} curvature and volume renormalization of {AHE} metrics on 4-manifolds},
        date={2001},
        ISSN={1073-2780},
     journal={Math. Res. Lett.},
      volume={8},
      number={1-2},
       pages={171\ndash 188},
      review={\MR{1825268}},
}

\bib{Avez1970}{article}{
      author={Avez, Andr\'{e}},
       title={Characteristic classes and {Weyl} tensor: {Applications} to general relativity},
        date={1970},
        ISSN={0027-8424},
     journal={Proc. Nat. Acad. Sci. U.S.A.},
      volume={66},
       pages={265\ndash 268},
         url={https://doi.org/10.1073/pnas.66.2.265},
      review={\MR{263098}},
}

\bib{BaileyEastwoodGraham1994}{article}{
      author={Bailey, Toby~N.},
      author={Eastwood, Michael~G.},
      author={Graham, C.~Robin},
       title={Invariant theory for conformal and {CR} geometry},
        date={1994},
        ISSN={0003-486X},
     journal={Ann. of Math. (2)},
      volume={139},
      number={3},
       pages={491\ndash 552},
         url={https://doi.org/10.2307/2118571},
      review={\MR{1283869}},
}

\bib{CaseLinYuan2022or}{article}{
      author={Case, Jeffrey~S},
      author={Lin, Yueh-Ju},
      author={Yuan, Wei},
       title={Curved versions of the {Ovsienko}--{Redou} operators},
        date={2023},
        ISSN={1073-7928,1687-0247},
     journal={Int. Math. Res. Not. IMRN},
      number={19},
       pages={16904\ndash 16929},
         url={https://doi.org/10.1093/imrn/rnad053},
      review={\MR{4651902}},
}

\bib{CaseYan2024}{article}{
      author={Case, Jeffrey~S.},
      author={Yan, Zetian},
       title={Formal self-adjointness of a family of conformally invariant bidifferential operators},
        date={2026},
        ISSN={0944-2669,1432-0835},
     journal={Calc. Var. Partial Differential Equations},
      volume={65},
      number={3},
       pages={Paper No. 100},
         url={https://doi.org/10.1007/s00526-026-03286-5},
      review={\MR{5033066}},
}

\bib{ChangGe2018}{article}{
      author={Chang, Sun-Yung~A.},
      author={Ge, Yuxin},
       title={Compactness of conformally compact {Einstein} manifolds in dimension 4},
        date={2018},
        ISSN={0001-8708,1090-2082},
     journal={Adv. Math.},
      volume={340},
       pages={588\ndash 652},
         url={https://doi.org/10.1016/j.aim.2018.10.010},
      review={\MR{3886176}},
}

\bib{ChangGe2020}{article}{
      author={Chang, Sun-Yung~A.},
      author={Ge, Yuxin},
      author={Qing, Jie},
       title={Compactness of conformally compact {Einstein} 4-manifolds {II}},
        date={2020},
        ISSN={0001-8708,1090-2082},
     journal={Adv. Math.},
      volume={373},
       pages={107325, 33},
         url={https://doi.org/10.1016/j.aim.2020.107325},
      review={\MR{4130465}},
}

\bib{ChangGeQing2020}{article}{
      author={Chang, Sun-Yung~A.},
      author={Ge, Yuxin},
      author={Qing, Jie},
       title={Compactness of conformally compact {Einstein} 4-manifolds {II}},
        date={2020},
        ISSN={0001-8708,1090-2082},
     journal={Adv. Math.},
      volume={373},
       pages={107325, 33},
         url={https://doi.org/10.1016/j.aim.2020.107325},
      review={\MR{4130465}},
}

\bib{ChangQingYang2006}{article}{
      author={Chang, Sun-Yung~A.},
      author={Qing, Jie},
      author={Yang, Paul.},
       title={On the renormalized volumes for conformally compact {Einstein} manifolds},
        date={2008},
     journal={J.\ Math.\ Sciences (N.\ Y.)},
      volume={149},
      number={6},
       pages={1755\ndash 1769},
      review={\MR{2336463}},
}

\bib{ChenLu2024}{article}{
      author={Chen, Fei-Yu},
      author={L\"u, H.},
       title={Holographic {Weyl} anomaly in {8D} from general higher curvature gravity},
        date={2025},
        ISSN={2470-0010,2470-0029},
     journal={Phys. Rev. D},
      volume={111},
      number={8},
       pages={Paper No. 086032, 16},
         url={https://doi.org/10.1103/physrevd.111.086032},
      review={\MR{4912725}},
}

\bib{Chern1955}{article}{
      author={Chern, Shiing-Shen},
       title={On curvature and characteristic classes of a {Riemann} manifold},
        date={1955},
        ISSN={0025-5858},
     journal={Abh. Math. Sem. Univ. Hamburg},
      volume={20},
       pages={117\ndash 126},
         url={https://doi.org/10.1007/BF02960745},
      review={\MR{0075647}},
}

\bib{ChruscielDelayLeeSkinner2005}{article}{
      author={Chru{\'s}ciel, Piotr~T.},
      author={Delay, Erwann},
      author={Lee, John~M.},
      author={Skinner, Dale~N.},
       title={Boundary regularity of conformally compact {Einstein} metrics},
        date={2005},
        ISSN={0022-040X},
     journal={J. Differential Geom.},
      volume={69},
      number={1},
       pages={111\ndash 136},
         url={http://projecteuclid.org/euclid.jdg/1121540341},
      review={\MR{2169584}},
}

\bib{FeffermanGraham2012}{book}{
      author={Fefferman, Charles},
      author={Graham, C.~Robin},
       title={The ambient metric},
      series={Annals of Mathematics Studies},
   publisher={Princeton University Press},
     address={Princeton, NJ},
        date={2012},
      volume={178},
        ISBN={978-0-691-15313-1},
      review={\MR{2858236}},
}

\bib{FullingKingWybourneCummins1992}{article}{
      author={Fulling, S.~A.},
      author={King, R.~C.},
      author={Wybourne, B.~G.},
      author={Cummins, C.~J.},
       title={Normal forms for tensor polynomials. {I}. {The} {Riemann} tensor},
        date={1992},
        ISSN={0264-9381},
     journal={Classical Quantum Gravity},
      volume={9},
      number={5},
       pages={1151\ndash 1197},
         url={http://stacks.iop.org/0264-9381/9/1151},
      review={\MR{1163881}},
}

\bib{Gilkey1975}{article}{
      author={Gilkey, Peter~B.},
       title={Local invariants of a pseudo-{Riemannian} manifold},
        date={1975},
        ISSN={0025-5521},
     journal={Math. Scand.},
      volume={36},
       pages={109\ndash 130},
         url={https://doi.org/10.7146/math.scand.a-11566},
      review={\MR{375340}},
}

\bib{GonzalezQing2010}{article}{
      author={Gonz{\'a}lez, Mar{\'{\i}}a del~Mar},
      author={Qing, Jie},
       title={Fractional conformal {Laplacians} and fractional {Yamabe} problems},
        date={2013},
        ISSN={2157-5045},
     journal={Anal. PDE},
      volume={6},
      number={7},
       pages={1535\ndash 1576},
         url={http://dx.doi.org/10.2140/apde.2013.6.1535},
      review={\MR{3148060}},
}

\bib{Graham2000}{inproceedings}{
      author={Graham, C.~Robin},
       title={Volume and area renormalizations for conformally compact {Einstein} metrics},
        date={2000},
   booktitle={The {P}roceedings of the 19th {W}inter {S}chool ``{G}eometry and {P}hysics'' ({S}rn\'\i, 1999)},
       pages={31\ndash 42},
      review={\MR{MR1758076}},
}

\bib{GJMS1992}{article}{
      author={Graham, C.~Robin},
      author={Jenne, Ralph},
      author={Mason, Lionel~J.},
      author={Sparling, George A.~J.},
       title={Conformally invariant powers of the {Laplacian}. {I}. {Existence}},
        date={1992},
        ISSN={0024-6107},
     journal={J. London Math. Soc. (2)},
      volume={46},
      number={3},
       pages={557\ndash 565},
      review={\MR{MR1190438}},
}

\bib{GrahamZworski2003}{article}{
      author={Graham, C.~Robin},
      author={Zworski, Maciej},
       title={Scattering matrix in conformal geometry},
        date={2003},
        ISSN={0020-9910},
     journal={Invent. Math.},
      volume={152},
      number={1},
       pages={89\ndash 118},
         url={http://dx.doi.org/10.1007/s00222-002-0268-1},
      review={\MR{1965361}},
}

\bib{Juhl2013}{article}{
      author={Juhl, Andreas},
       title={Explicit formulas for {GJMS}-operators and {$Q$}-curvatures},
        date={2013},
        ISSN={1016-443X},
     journal={Geom. Funct. Anal.},
      volume={23},
      number={4},
       pages={1278\ndash 1370},
         url={http://dx.doi.org/10.1007/s00039-013-0232-9},
      review={\MR{3077914}},
}

\bib{Maldacena1998}{article}{
      author={Maldacena, Juan},
       title={The large {$N$} limit of superconformal field theories and supergravity},
        date={1998},
        ISSN={1095-0761},
     journal={Adv. Theor. Math. Phys.},
      volume={2},
      number={2},
       pages={231\ndash 252},
      review={\MR{MR1633016}},
}

\bib{Marugame2017}{article}{
      author={Marugame, Taiji},
       title={Some remarks on the total {CR} {$Q$} and {$Q'$}-curvatures},
        date={2018},
        ISSN={1815-0659},
     journal={SIGMA Symmetry Integrability Geom. Methods Appl.},
      volume={14},
       pages={Paper No. 010, 8},
         url={https://doi.org/10.3842/SIGMA.2018.010},
      review={\MR{3763373}},
}

\bib{Witten1998}{article}{
      author={Witten, Edward},
       title={Anti de {Sitter} space and holography},
        date={1998},
        ISSN={1095-0761},
     journal={Adv. Theor. Math. Phys.},
      volume={2},
      number={2},
       pages={253\ndash 291},
      review={\MR{MR1633012}},
}

\end{biblist}
\end{bibdiv}
\end{document}